\documentclass[12pt]{amsart}
\newenvironment{prf}{\noindent{\it Proof\/}:}{$\;\square$
\par\medskip}

\newtheorem%
{thm}{Theorem}[section]
\newtheorem
{proposition}[thm]{Proposition}
\newtheorem%
{lemma}[thm]{Lemma}
\newtheorem%
{lemmadef}[thm]{Lemma-Definition}
\newtheorem%
{corollary}[thm]{Corollary}
\newtheorem%
{conjecture}[thm]{Conjecture}

\newcommand{\dontprint}[1]{\relax}
\newcommand{\ndash}{\nobreakdash-\hspace{0pt}}
\newcommand{\Ndash}{\nobreakdash--\hspace{0pt}}
\newcommand{\eps}{\varepsilon}
\title
[Relative formality theorem]
{Relative formality theorem and quantisation of coisotropic submanifolds}
\author{Alberto S. Cattaneo and Giovanni Felder}
\address{Institut f\"ur Mathematik, Universit\"at Z\"urich--Irchel,  
Winterthurerstrasse 190, CH-8057 Z\"urich, Switzerland}  
\email{alberto.cattaneo@math.unizh.ch}
\address{D-MATH, ETH-Zentrum, CH-8092 Z\"urich, Switzerland}
\email{giovanni.felder@math.ethz.ch}

\thanks{A.~S.~C. acknowledges partial support of SNF Grant No.~20-100029/1}
\thanks{G. F. acknowledges partial support of SNF Grant No.~200020-105450/1}
                                                                               
\begin{document}

\begin{abstract} 
We prove a relative version of Kontsevich's formality theorem. This
theorem involves a manifold $M$ and a submanifold $C$ and reduces to
Kontsevich's theorem if $C=M$. It states that the DGLA of
multivector fields on an infinitesimal neighbourhood of $C$ is
$L_\infty$\ndash quasiisomorphic to the DGLA of multidifferential
operators acting on sections of the exterior algebra of the conormal
bundle.  Applications to the deformation quantisation of coisotropic
submanifolds are given. The proof uses a duality transformation to
reduce the theorem to a version of Kontsevich's theorem for
supermanifolds, which we also discuss.  In physical language, the
result states that there is a duality between the Poisson
sigma model on a manifold with a D-brane and the Poisson sigma
model on a supermanifold without branes (or, more properly, with a
brane which extends over the whole supermanifold).
\end{abstract}
\dontprint{
We prove a relative version of Kontsevich's formality theorem. This
theorem involves a manifold M and a submanifold C and reduces to
Kontsevich's theorem if C=M. It states that the DGLA of
multivector fields on an infinitesimal neighbourhood of C is
L-infinity-quasiisomorphic to the DGLA of multidifferential
operators acting on sections of the exterior algebra of the conormal
bundle.  Applications to the deformation quantisation of coisotropic
submanifolds are given. The proof uses a duality transformation to
reduce the theorem to a version of Kontsevich's theorem for
supermanifolds, which we also discuss.  In physical language, the
result states that there is a duality between the Poisson
sigma model on a manifold with a D-brane and the Poisson sigma
model on a supermanifold without branes (or, more properly, with a
brane which extends over the whole supermanifold).
}

\maketitle
\section{Introduction}
In \cite{Kontsevich} Kontsevich gave a solution to the problem of
deformation quantisation of the algebra of functions on an arbitrary  Poisson manifold.
This solution is based on his formality theorem, stating that the differential
graded Lie algebra (DGLA) of multidifferential operators is $L_\infty$\ndash quasiisomorphic
to its cohomology, the DGLA of multivector fields. We consider here a version of
the formality theorem for a pair $(M,C)$ of manifolds $C\subset M$, which reduces to
the original formality theorem if $C=M$. The algebra of functions on $M$ is replaced
here by the graded commutative algebra $A$ of sections of the exterior algebra of the
normal bundle $NC$. The (suitably completed)
 Hochschild complex of $A$, with Hochschild differential
and Gerstenhaber bracket contains the sub\ndash DGLA $\hat{\mathcal{D}}(A)$ of ``multidifferential
operators'' on $A$, namely cochains built out of products of compositions of derivations of $A$. 
The statement is that this DGLA is $L_\infty$\ndash quasiisomorphic to its
cohomology, which is identified with the DGLA $\mathcal{T}(M,C)$ of multivector fields
on a formal neighbourhood of $C$ with Schouten\Ndash Nijenhuis bracket and zero
differential, see Theorem \ref{t-rft}. The proof is based on a ``Fourier transform'' 
Theorem \ref{t-Fields} which states that the Gerstenhaber algebra  $\mathcal{T}(M,C)$
is isomorphic to the Gerstenhaber algebra of multivector fields on the supermanifold
$N^*[1]C$, the conormal bundle with shifted parity of the fibres. In terms of supermanifolds, this
isomorphism is obtained from an isomorphism
of odd\ndash symplectic graded supermanifolds 
$T^*[1]NC\to T^*[1]N^*[1]C$, a variant of
what Roytenberg calls ``Legendre transform'' \cite{Roytenberg}. 
At this point the
result follows from a version of Kontsevich's formality theorem for supermanifolds,
see Theorem \ref{t-Holmes}. The proof of the latter theorem is parallel to Kontsevich's \cite{Kontsevich},
except for signs, which are already non\ndash trivial in the original setting of ordinary manifolds.
For this reason we work out all signs in the Appendix, and develop a formalism in which
these signs appear in an essentially transparent way.   

In the application to deformation quantisation we take $C$ to be a coisotropic
submanifold of $M$, which means
 that the vanishing ideal $I(C)$ of $C$ is a Poisson subalgebra
of $C^\infty(M)$. To these data one associates a Poisson algebra, the algebra
of functions on the reduced phase space $C^\infty(\underline C)=N(I(C))/I(C)$, the quotient of the normaliser
of the Lie algebra $I(C)$ by $I(C)$. Even if the reduced phase space $\underline C$, which is by
definition the space of leaves of the characteristic foliation of $C$, is singular,
the Poisson algebra  $C^\infty(\underline C)$ of ``smooth functions'' on it 
makes sense, and one can ask the
question of quantising this algebra in the sense of deformation quantisation. 
It seems that this is not always possible
because of anomalies, but it may be argued that the question is not correct and that
one should not try to quantise $C^\infty(\underline C)$, which can be trivial even
for interesting $C$, but rather some kind of resolution of it. In fact there is
a natural complex whose cohomology in degree zero is $C^\infty(\underline C)$:
the conormal bundle $N^*C$ of a coisotropic manifold is naturally a Lie algebroid,
and $C^\infty(\underline C)$ is its zeroth Lie algebroid cohomology algebra. Thus
one replaces $C^\infty(\underline C)$ by the Lie algebroid cochain complex
$\Gamma(C,\wedge NC)$. This differential graded algebra is however not a
Poisson algebra; it turns out, as essentially noticed by Oh and Park \cite{OP},
 that the Poisson structure on $M$ induces
a $P_\infty$\ndash structure on $\Gamma(C,\wedge NC)$, namely an $L_\infty$\ndash
structure whose structure maps are multiderivations, see Theorem \ref{t-Lowell}. 
Algebraically, the $P_\infty$\ndash brackets are obtained from the Poisson structure
on $M$ as higher derived brackets in the sense of T. Voronov \cite{Voronov}.
This $L_\infty$\ndash
structure induces the Poisson bracket on cohomology. At this point the
$L_\infty$\ndash machinery can be applied: the $L_\infty$\ndash structure can
be understood as a solution of the Maurer\Ndash Cartan equation in $\mathcal{T}(M,C)$
and is mapped by the $L_\infty$\ndash quasiisomorphism to a solution of the
Maurer\Ndash Cartan equation in $\hat{\mathcal{D}}(A)$. The latter is a deformation
of the product in $A$ as an $A_\infty$\ndash algebra, see Theorem \ref{t-azeroinfty}. 
At this point one may want to
pass to cohomology to quantise the original Poisson algebra $C^\infty(\underline C)$,
or more general the whole Lie algebroid cohomology. It is here that one meets
the anomalies in general. Namely the $A_\infty$\ndash algebra obtained by this
construction is not flat in general, namely its $0$th product map $\mu_0$ may not vanish
(we use the not-quite-standard but natural definition of $A_\infty$\ndash algebra
which allows for non-zero product $\mu_0\in A$ and call an $A_\infty$\ndash
algebra flat if $\mu_0=0$). In this case the first product $\mu_1$ is not
a differential and the cohomology is not defined. Removing $\mu_0$ (which is at least quadratic in the deformation parameter)
is a
cohomological problem with an obstruction in the second Lie algebroid cohomology
group. In some cases the obstruction vanishes even if the cohomology
does not, see the second remark in \ref{s-qua}. If the obstruction vanishes,
one gets an associative algebra, which is however not always a
flat deformation of $C^\infty(\underline C)$. This time the obstruction
is in the first cohomology group, see Corollary \ref{c-Longfellow}.

Some of the results presented here were announced in \cite{CaFe2003}. There the
interpretation of these results in terms of topological quantum field
theory is given: the $L_\infty$\ndash quasiisomorphism is constructed
using a topological sigma model on the disk with the boundary condition
that the boundary is sent to $C$. 
An alternative approach to this class of problems, based on Tamarkin's formality theorem, 
was proposed recently in \cite{BGHHW}. See also the very recent
preprint \cite{LS}, in which a use of the formality for supermanifolds
similar to ours is presented in a physics context and shown to be applicable
to weak Poisson manifolds.

{}From the point of view of topological quantum field theory adopted in
\cite{CaFe2003} this paper concerns the case of a single D-brane. The more general
case of several D-branes will be studied elsewhere. It corresponds to the
theory of (bi)modules over the deformed algebras. 
\medskip

\noindent{\bf Acknowledgements} We
are grateful to Jim Stasheff and Riccardo Longoni for useful 
comments and
corrections to the the first draft of this paper. 
We thank Martin Bordemann for his useful comments and explanations.
We also thank the referee for carefully reading the manuscript and
for suggesting improvements.
\medskip

\noindent{\bf Conventions.}
We work in the category of graded vector spaces (or free modules over a commutative
ring) $V=\oplus_{j\in\mathrm{Z}} V^j$,
and denote by $|a|$ the degree of a homogeneous element $a\in V^{|a|}$. We denote by $V[n]$ the
graded vector space $\oplus_j V[n]^j$ with 
$V[n]^j=V^{n+j}$. The
space of homomorphisms $f\colon V\to W$ of degree $j$ (i.e., such that $f(V^i)\subset W^{j+i}$ is denoted by $\mathrm{Hom}^j(V,W)$. The Koszul
sign rule holds. A derivation $f$ of degree $|f|$ of a graded algebra $A$ 
is a linear endomorphisms of degree $|f|$ 
obeying $f(ab)=f(a)b+(-1)^{|a||f|}af(b)$ for all $a,b\in A$. 
See \ref{N&C} for more details.

\section{Coisotropic submanifolds of Poisson manifolds}
\subsection{Coisotropic submanifolds}
Let $(M,\pi)$ be a Poisson manifold,
with Poisson bivector field $\pi\in\Gamma(M,\wedge^2TM)$ and
Poisson bracket $\{f,g\}=\langle\pi,df\otimes dg\rangle$.
Let $\pi^\sharp\colon T^*M\to TM$ be the bundle map induced by $\pi$
on each cotangent space: $\langle\pi^\sharp(\alpha),\beta\rangle
=\langle\pi,\alpha\otimes\beta\rangle$. A submanifold
$C\subset M$ is called {\em coisotropic} \cite{Weinstein}
if $\pi^\sharp|_C$
maps the conormal bundle $N^*C=\mathrm{Ann}(TC)\subset T^*_CM$ 
to the tangent bundle $TC$. 
Equivalently, $C$ is coisotropic if and only if the
ideal $I(C)$ of the algebra $C^\infty(M)$ 
consisting of functions vanishing on $C$ is closed
under the Poisson bracket.
Examples include $M$ itself, Lagrangian submanifolds of
symplectic manifolds, graphs of Poisson maps,
zeros of equivariant moment maps
 and 
mechanical systems with first class constraints.

Coisotropic submanifolds come with interesting geometric
and algebraic structures, which we turn to describe.

\subsection{Characteristic foliation and reduced phase space} \label{s-HR}
If $C$ is coisotropic, the distribution $\pi^\sharp(N^*C)\subset
TC$ of tangent subspaces is involutive since it is spanned
by hamiltonian vector fields $X_h=\pi^\sharp dh$ with $h\in I(C)$,
which commute on $C$ by the coisotropy condition.
The corresponding
foliation is the {\em characteristic foliation} of $C$. The leaves of the characteristic foliation have points
related by hamiltonian flows
with hamiltonian functions in $I(C)$.
The {\em reduced phase space} $\underline C$ of 
$C$ is the space of leaves of the characteristic foliation.
It can be a wild space so that it is better to consider the
algebra of functions on it, which is, by definition, the
algebra of functions on $C$ that are invariant under the hamiltonian flows 
of $I(C)$:
\[
C^\infty(\underline C)=\{f\in C^\infty(C)\,|\,X_h(f)=0\quad
\forall h\in I(C)\}.
\]
If $\underline C$ is a manifold, then the Poisson bivector
field descends to $\underline C$ and gives it a structure
of Poisson manifold. In general $C^\infty(\underline C)$
is a Poisson algebra, namely a commutative algebra with
a Lie bracket which is a derivation in each of its arguments. 
The Lie algebra structure is induced
from the Lie algebra structure on $C^\infty(M)$
as is clear from the representation
\[
C^\infty(\underline C)=N(I(C))/I(C),
\]
as the quotient by $I(C)$ of the normaliser 
\[N(I(C))=\{f\in C^\infty(M)\,|\,
\{I(C),f\}\subset I(C)\}\]
of the Lie subalgebra $I(C)$.

\subsection{The Lie algebroid of a coisotropic submanifold}
The space of sections of the conormal bundle 
has a natural Lie algebra structure. The Lie bracket 
is uniquely defined
by the conditions
\begin{eqnarray*}
[df,dg]&=&d\{f,g\},\qquad f, g\in I(C),
\\
{}[f\alpha,\beta]&=&f[\alpha,\beta]-
\pi^\sharp(\beta)(f)\alpha,\qquad 
\alpha,\beta\in \Gamma(C,N^*C),\quad f\in C^\infty(C).
\end{eqnarray*}
By construction, $\pi^\sharp$ induces a Lie algebra homomorphism
$\Gamma(C,N^*C)\to \Gamma(C,TC)$ from this Lie algebra to
the Lie algebra of vector fields. In other words, $N^*C$ is
a Lie algebroid over $C$.
\subsection{The cochain complex of a coisotropic submanifold}\label{s-cochain}
As for every Lie algebroid, the Lie algebroid of a coisotropic
submanifold comes with a cochain complex, see \cite{Mackenzie}:
\[
\cdots\to
\Gamma(C,\wedge^jNC)\to\Gamma(C,\wedge^{j+1}NC)\to\cdots.
\]
The differential $\delta$ on $\Gamma(C,\wedge^0NC)=C^\infty(C)$ is $\delta f=\pi^\sharp d\tilde f \mod TC$, for any extension
$\tilde f$ of $f$ to $M$: the class of $\pi^\sharp d\tilde f$
in $NC=T_CM/TC$ is independent of the choice of extension because
of the coisotropy condition. The differential on $\Gamma(C,
\wedge^1NC)$ is the dual map to the Lie bracket $\Gamma(C,\wedge^2N^*C)\to\Gamma(C,N^*C)$. The differential on general cochains
is determined by the rule
\[
\delta(\alpha\wedge\beta)=\delta\alpha\wedge\beta
+(-1)^{|\alpha|}\alpha\wedge\delta\beta,\qquad\alpha,\beta\in
\Gamma(C,\wedge NC).
\]
The cohomology of this complex is the cohomology
$H_\pi(N^*C)$ of the Lie algebroid $N^*C$. It is a graded
commutative algebra. In degree $0$ we
have
\[
H^0_\pi(N^*C)=C^\infty(\underline C).
\]
The first cohomology group describes infinitesimal deformations
of the imbedding of $C$ as a 
coisotropic submanifold up to deformations induced by
hamiltonian flows.

\subsection{The $P_\infty$\ndash structure on the cochain complex}\label{s-PINF}
A natural question is whether the Poisson bracket
on $H^0_\pi(N^*C)=C^\infty(\underline C)$ comes
from a structure on the cochain complex.
We want to show that this structure is
a flat $P_\infty$\ndash structure (defined up to homotopy),
namely a graded commutative algebra structure with
a compatible flat $L_\infty$\ndash structure. In particular,
this flat $P_\infty$\ndash structure induces a Poisson bracket
on the whole cohomology algebra $H_\pi(N^*M)$.
 The definitions
are as follows.
A \dontprint{weak}
$P_{\infty}$\ndash {\em algebra} ($P$ for Poisson) is a graded
commutative algebra $A$ over a field of characteristic zero
with a sequence of linear
maps $\lambda_n\colon A^{\otimes n}\to A$ of degree $2-n$, 
$n=0,1,2,\dots$ with the following properties (Properties (i) and
(iii) characterise $L_\infty$\ndash algebras) that are to hold for
arbitrary $a_1,\dots,a_n\in A$:

\medskip

\noindent{(i)}
$\lambda_n(\dots,a_i,a_{i+1},\dots)=
-(-1)^{|a_i|\cdot|a_{i+1}|}\lambda_n(\dots,a_{i+1},a_i,\dots).$

\noindent{(ii)}
$a\mapsto\lambda_n(a_1,\dots,a_{n-1},a)$ is a derivation of degree
$2-n-\sum_{i=1}^{n-1} |a_i|$.

\noindent{(iii)}
For all $n\geq 0$, the map
\[
a_1\otimes\cdots\otimes a_n\mapsto
\sum_{q=0}^{n} \frac{(-1)^{q(n-q)}}{q!(n-q)!}\lambda_{n-q+1}(\lambda_q(a_1,\dots,a_q),
a_{q+1},\dots,a_n)
\]
vanishes on the image of the alternation map $\mathrm{Alt}_n
=\frac1{n!}\sum_{\sigma\in S_n}\mathrm{sign}(\sigma)\sigma$.

\medskip

A {\em flat} $P_\infty$\ndash algebra is a $P_{\infty}$\ndash algebra such that $\lambda_0=0$.
Then $\lambda_1$ is a differential, $\lambda_2$ a chain map obeying the
Jacobi identity up to exact terms. So the $\lambda_1$\ndash cohomology of a
flat $P_\infty$\ndash algebra is a graded Poisson algebra.

\medskip\noindent{\bf Remark.} This notion of flat $P_\infty$-structure is not completely
standard: in the spirit of homotopical algebra one might prefer a notion in which also
the commutativity of the product and the Leibniz rule hold only up to homotopy. Here
only the bracket of a Poisson algebra is replaced by a sequence of higher brackets
controlling the violation of the Jacobi identity.

\medskip

The construction of a flat $P_\infty$\ndash structure on $A=\Gamma(C,\wedge N C)$ 
for a coisotropic submanifold $C\subset M$ depends
on the choice of an identification of a tubular neighbourhood of $C$ with
the normal bundle of $C$, more precisely an embedding $\iota$ of $NC$ into $M$
sending the zero section identically
to $C$ and such that, for $x\in C$, the restriction of
$\iota_*\colon T_x(NC)\to T_xM$ 
composed with the canonical projection
is the identity $N_xC\to T_xM/T_xC$. 
As such an embedding is unique
up to homotopy, the construction gives a flat $P_\infty$\ndash structure up to
homotopy. We give the construction in the more general setting of
a general submanifold of $M$, yielding a non\ndash necessarily
flat $P_{\infty}$\ndash structure.
As the construction only involves a neighbourhood of the
submanifold, we may as well assume that $M$ is the total
space of a vector bundle, which is then canonically
the normal bundle of its zero section $C$.

\begin{proposition}\label{p-Pzeroinfty}
 Let $C$ be a submanifold of a Poisson manifold $M$, not
necessarily coisotropic.
Assume that $M$ is the total space of a vector bundle $p\colon E\to C$ and let $C\subset M$ 
be the zero section of $E$. Then there is
a unique $P_{\infty}$\ndash structure on $\Gamma(C,\wedge NC)\simeq\Gamma(C,\wedge E)$ 
such that for $v_1,\dots,v_n\in\Gamma(C,E)$, 
 $f,g\in C^\infty(C)=\Gamma(C,\wedge^0E)$, and
$u,w\in\Gamma(C,E^*)$
\begin{eqnarray*}
\lambda_n(v_1,\dots,v_{n-2},f,g)&=&(-1)^{n-2}v_1\cdots v_{n-2}\{p^*f,p^*g\}|_C,\\
\langle\lambda_n(v_1,\dots,v_{n-1},f),u\rangle
&=&(-1)^{n-1}v_1\cdots v_{n-1}\{p^*f,u\}|_C\\
\langle\lambda_n(v_1,\dots,v_{n}),u\otimes w\rangle
&=&(-1)^{n}v_1\cdots v_{n}\{u,w\}|_C.
\end{eqnarray*}
On the right-hand side of these equations, $u$, $w$ are regarded as 
functions on $M$ linear on the fibres, and $v_i$ as vertical vector fields on $M$.
\end{proposition}

The uniqueness part is clear: since $\lambda_n$ is a multiderivation, it is
sufficient to define it on $\Gamma(C,\wedge^jE)$ with $j=0,1$ as the algebra
is generated by these spaces. Since $\lambda_n$ is of degree $2-n$, $\lambda_n$
vanishes on elements of degree 0 or 1 except in the three cases listed in
the Proposition. The fact that the $\lambda_n$ extend to a 
\dontprint{weak} $P_{\infty}$\ndash structure
can be checked directly, but we will deduce it from a more general result below
(see Prop.~\ref{p-Dalhousie}). 

\begin{thm}\label{t-Lowell}
Assume that in Prop.~\ref{p-Pzeroinfty}, $C$ is coisotropic.
Then $\lambda_0=0$ so that $\Gamma(C,\wedge NC)$ is
a flat $P_\infty$\ndash algebra, $\lambda_1$ is the differential of \ref{s-cochain}
and $\lambda_2$ induces the Poisson bracket on $H^0_\pi(C)=C^\infty(\underbar C)$
of \ref{s-HR}.
\end{thm}

\begin{prf} Let $u,v\in\Gamma(C,N^*C)$ considered as fibre\ndash linear functions on $NC$.
We have $\langle\lambda_0,u\otimes v\rangle=\{u,v\}|_C=0$ since $u,v$ belong to
the ideal $I(C)$, which is a Lie subalgebra if $C$ is coisotropic. By definition,
if $f\in C^\infty(C)$,
$\langle\lambda_1(f),u\rangle=\{p^*f,u\}|_C$, so $\lambda_1(f)$ is the class
in $\Gamma(C,NC)$ of $\pi^\sharp d\tilde f$ with $\tilde f=p^*f$. If $w\in\Gamma(C,NC)$, 
   $\langle\lambda_1(w),u\otimes v\rangle=w\{u,v\}|_C=\langle w,[u,v]\rangle$. Therefore
$\lambda_1$ is indeed the differential of \ref{s-cochain}. 
As for $\lambda_2$, recall that the Poisson bracket of $f,g\in 
H^0_\pi(N^*C)=\mathrm{Ker}(\lambda_1\colon C^\infty(C)\to \Gamma(C,E))$ is
defined as $\{\tilde f,\tilde g\}$ for any extension $\tilde f,\tilde g$ of
$f,g$ to $M$. This is precisely the definition of $\lambda_2$, with $\tilde f=p^*f$.
\end{prf}

Using this result, we obtain a Poisson bracket  induced by $\lambda_2$ on the cohomology
$H_\pi(N^*C)$.  A priori this bracket depends on
the choice of embedding of $NC$ into $M$.
However we see by a standard homotopy argument that this is
not the case:

\begin{proposition} Let $C\subset M$ be coisotropic. Then the Lie
bracket induced by $\lambda_2$ on the cohomology $H_\pi(N^*C)$ is
independent of the choice of embedding of $NC$ 
into a tubular neighbourhood
of $C$
\end{proposition}

\subsection{Higher derived brackets and relative multivector fields}
The maps $\lambda_j$ on $\Gamma(C,\wedge NC)$
are a special case of higher derived brackets, see
\cite{Voronov}. Let $\mathfrak a$ be an abelian graded Lie subalgebra of a graded 
Lie algebra $\mathfrak g$, with a projection $P\colon\mathfrak g\to \mathfrak a$,
satisfying $P[a,b]=P[Pa,b]+P[a,Pb]$. Suppose we have an
element $\pi\in \mathfrak g$ of degree $1$ obeying $[\pi,\pi]=0$. Then the higher derived
brackets
\[
\{a_1,\dots,a_n\}=P[\cdots[\pi,a_1],a_2],\dots],a_n]
\]
are graded {\em symmetric} multilinear functions of degree 1, obeying
the Jacobi identities
\[
\sum_k
\sum_{\sigma\in S_n}\frac{(-1)^\epsilon}{k!(n-k)!}\{\{a_{\sigma(1)},\dots,a_{\sigma(k)}\},\dots,a_{\sigma(n)}\}, 
\]
with the natural sign: $\epsilon=\sum_{i<j,\sigma(i)>\sigma(j)}\mathrm{deg}(a_i)
\mathrm{deg}(a_j)$.
The brackets
\[
\lambda_n(a_1,\dots,a_n)=(-1)^{\sum_i (i-1)\mathrm{deg}(a_i)}\{a_1,\dots,a_n\} 
\]
are then {\em skew-symmetric} and of degree $2-n$ 
with respect to the shifted degree $|a|=\mathrm{deg}(a)+1$
and they obey the $L_\infty$\ndash Jacobi identities (iii) above.

In our case, $\mathfrak g=\mathcal{T}(M,C)$ is the Lie algebra of {\em relative multivector fields} on the submanifold $C\subset M$ with the Schouten\Ndash Nijenhuis
bracket. It is the inverse limit $\varprojlim\mathcal{T}(M)/I(C)^n\mathcal{T}(M)$ 
where $\mathcal{T}(M)=\oplus_{j=-1}^{\infty}\mathcal{T}^j(M)$ is the graded Lie algebra
of multivector fields and $I(C)$ is the ideal in $C^\infty(M)$ of functions vanishing
on $C$. The Lie subalgebra $\mathfrak a=\Gamma(C,\wedge NC)$ consists of sums of products of
vector fields tangent to the fibres of $M=E\to C$ and constant along
each fibre.

\section{Quantisation of coisotropic submanifolds}

\subsection{$A_{\infty}$\ndash algebras and flat $A_\infty$\ndash algebras} 
An $A_{\infty}$\ndash algebra \cite{Stasheff63}
over a commutative ring $R$ is a free graded left $R$\ndash module $A=\oplus_{j\in\mathbb Z}A^j$
with $R$-linear maps $\mu_n\colon A^{\otimes n}\to A[2-n]$ of degree $0$ ($n=0,1,\dots$) obeying
the associativity relations
\begin{eqnarray*}
\lefteqn{\sum_{q=0}^{n} (-1)^{q(n-q)}\sum_{j=0}^{p-1} (-1)^{(q-1)j+\sum_{i=1}^j|a_i|\cdot q}}\\
&&\mu_{n-q+1}(a_1,\dots,a_j,\mu_q(a_{j+1},\dots,a_{j+q}),
a_{j+q+1},\dots,a_n)=0.
\end{eqnarray*}
A {\em flat} $A_\infty$\ndash algebra is an $A_{\infty}$\ndash algebra with $\mu_0=0$. In this case
$\mu_1$ is a differential and $\mu_2$ induces an associative product on its 
cohomology. Associative algebras can be regarded as $A_\infty$\ndash algebras with
product $\mu_2$ and all other $\mu_i=0$. 

Let $F(V)=\oplus_{j\geq 0} V^{\otimes j}$ be the tensor coalgebra over $R$ 
generated by a free $R$-module $V$ and denote by $p_j:F(V)\to V^{\otimes j}$ 
the projection onto the $j$th summand. Let $A[1]$ be the graded $R$-module
with homogeneous components $A[1]^j=A^{j+1}$ and let $s:A[1]\to A$ be the tautological
map of degree 1.
Then an $A_\infty$-structure on $A$ is the same as a coderivation 
$Q$ of degree 1 of $F(A[1])$ obeying $[Q,Q]=0$, see \ref{s-Hoch} in the Appendix. 
The coderivation $Q$ and
the products $\mu_n$ are related by 
\[
\mu_n\circ(s\otimes\cdots\otimes s)=s\circ p_1\circ Q|_{A[1]^{\otimes n}}.
\]
The ``strange'' signs in the associativity relations come from the Koszul
rule, if we take into account that $s$ has degree 1.

The following result is a graded, $A_\infty$\ndash 
version of the classical result relating
first order associative deformations of algebras of smooth functions to Poisson brackets. 

\begin{proposition}
Let $A_0=\oplus_i\Gamma(M,\wedge^i E)$  be the graded commutative
 algebra of sections of the
exterior algebra of a vector bundle $E\to M$. Let $(\mu_n)_{n=0}^\infty$
be an $A_\infty$\ndash algebra structure
on $A=A_0[[\epsilon]]$ over $\mathbb{R}[[\epsilon]]$ which reduces modulo $\epsilon$
to the algebra structure on $A_0=A/\epsilon A$ and
such that the structure maps $\mu_n$ are multidifferential operators.
Let $\lambda_n=\frac1\epsilon\mu_n\circ\sum_{\sigma\in
S_n}\mathrm{sign}(\sigma)\sigma
\mod \epsilon A$.
 Then $(\lambda_n)_{n=0}^\infty$ is a $P_\infty$\ndash
structure on $A_0$.
\end{proposition}

\begin{prf}
The products have the form
$\mu_n=\mu^{0}_n+\epsilon\mu^{1}_n+\cdots$ with $\mu_n^{0}=0$ except
for $n=2$. The associativity relations can be expressed as $[\mu,\mu]=0$ interms of the Gerstenhaber bracket (see \eqref{e-GB} below), so that to lower
order in $\epsilon$ we have
\[ 
[\mu^0,\mu^0]=0,\qquad [\mu^0,\mu^1]=0,\qquad [\mu^1,\mu^1]+2[\mu^0,\mu^2]=0.
\]
The first equation is just the associativity of the product in $A_0$.
The second equation states that $\mu_p^1$ is a Hochschild cocycle for
each $p$:
$b\mu^1_p=0$. By the HKR theorem (see Lemma
\ref{l-HKR}), we have 
$\mu^1_p=\tilde\mu^1_p+b\varphi_p$ for some
alternating multiderivation $\tilde\mu^1$. A straightforward
direct calculation shows that, for any $p$-cochain
$\varphi_p$,
\[
b\varphi_p\circ\mathrm{Alt}_{p+1}=0,
\]
owing to the commutativity of the product.
Thus $\lambda=\tilde\mu^1$ and therefore $\lambda$ is an alternating
multiderivation, i.e., obeys (i) (ii).
Finally, the third equation, restricted to skew-symmetric tensors in $A_0^{\otimes n}$
becomes $[\lambda,\lambda]=0$ as the term $[\mu^0,\mu^2]=b\mu^2$ does not contribute, again
because of the commutativity of $\mu^0$. This proves property (iii) of $P_\infty$-structures.
\end{prf}

\subsection{Quantisation}\label{s-qua}
Our problem is to quantise the Poisson algebra $H_\pi(N^*C)$ (or at
least the subalgebra $H^0_\pi(N^*C)=C^\infty(\underline C)$) 
for coisotropic $C$, namely to
find a {\em star\ndash product}, i.e., 
an associative $\mathbb R[[\epsilon]]$-bilinear product $\star$
on $H_\pi(N^*C)[[\epsilon]]$ deforming the graded commutative product and
such that $\epsilon^{-1}(a\star b-(-1)^{|a|\cdot|b|}b\star a)$ 
is the Poisson bracket modulo $\epsilon$.
It seems that this is impossible in general. What one always
has is a quantisation
of the $P_{\infty}$\ndash algebra $\Gamma(C,N^*C)$ for any submanifold $C$ as
an $A_{\infty}$\ndash algebra. From this result we then solve the original
quantisation problem if suitable obstructions vanish.

\begin{thm}\label{t-azeroinfty}
Let $C\subset M$ be a submanifold of a Poisson manifold $M$ and let
 $A=\Gamma(C,\wedge NC)$ with a $P_\infty$\ndash structure $\lambda$ induced
by the Poisson structure on $M$. Then there is an $\mathbb{R}[[\epsilon]]$-linear
$A_{\infty}$\ndash structure on $A_\epsilon=A[[\epsilon]]$ inducing $\lambda$
on $A_\epsilon/\epsilon A_\epsilon$. 
Moreover,
if $C$ is coisotropic (so that $\lambda_0=0$) 
then $\mu_0=O(\epsilon^2)$. 
\end{thm}

In the case where $C$ is an affine subspace in $\mathbb{R}^n$ there
is an explicit Feynman diagram expansion describing this $A_{\infty}$\ndash structure,
see \cite{CaFe2003}. 

Theorem \ref{t-azeroinfty} is proved below as
a consequence of the relative formality Theorem \ref{t-rft}.

\begin{corollary}\label{c-Longfellow} If $C$ is coisotropic and $H^2_\pi(N^*C)=0$ then the 
$A_{\infty}$\ndash structure in the preceding theorem may be chosen as 
a flat $A_\infty$\ndash structure $\mu=(\mu_1,\mu_2,\dots)$. In particular,
$\mu_2$ induces an associative product on the cohomology $H_{\pi,\epsilon}(N^*C)$
of the complex $(A[[\epsilon]],\mu_1)$.
If additionally $H^1_\pi(N^*C)=0$, then there is an isomorphism of 
$\mathbb{R}[[\epsilon]]$\ndash modules $H^0_{\pi,\epsilon}(N^*C)\to H^0_\pi(N^*C)[[\epsilon]]$
sending $\mu_2$ to a star-product on the Poisson algebra $H^0_\pi(N^*C)$.
\end{corollary}

\begin{prf}
Theorem \ref{t-azeroinfty} gives an $A_\infty$\ndash structure $(\mu_n)_{n\geq0}$
with $\mu_0=O(\epsilon^2)$ and $\mu_n=O(\epsilon)$ for $n\neq2$. 
The problem
is to find an $a\in\epsilon A[[\epsilon]]$ of degree 1 such that 
\begin{equation}\label{e-PietroBachi}
\sum_{n=0}^\infty\mu_n(a,a,\dots,a)=0.
\end{equation}
Suppose for the moment that we have such an $a$.
Let $Q$ be the coderivation of $F(A[1])$ associated to $\mu$ and denote by 
$T\in\mathrm{End}(F(A[1]))$ 
the unique coalgebra automorphism with vanishing 
Taylor components $T_j=p_1\circ T|_{A[1]^{\otimes j}}$ except
 $T_0=a$ and $T_1(x)=x$, $x\in A[1]$.
Then $\hat Q=T^{-1}\circ Q\circ T$ is a coderivation of $F(A[1))$ defining
a new $A_\infty$\ndash structure $\hat \mu$
obeying the properties of Theorem \ref{t-azeroinfty},
but with $\hat\mu_0=0$. Indeed, $\hat\mu_0=p_1\circ\hat Q|_{V^{\otimes 0}}
=\sum p_1\circ T_1^{-1}\circ 
Q\circ(T_0\otimes\cdots\otimes T_0)
=\sum_{n\geq0} Q(a,\dots,a)=0$.

Equation \eqref{e-PietroBachi} for a power series
$a=\epsilon a_1+\epsilon^2 a_2+\cdots\in\epsilon\Gamma(C,NC)[[\epsilon]]$ 
can then be solved in a standard
recursive way as a cohomological problem for the 
Lie algebroid differential $d=\epsilon^{-1}\mu_1|_{\epsilon=0}=\lambda_1$.
Since $a$ is odd and $\mu_2$ is commutative to lowest order, 
we have $\mu_2(a,a)=O(\epsilon^3)$, 
so if $\mu_0=\epsilon^2F+O(\epsilon^3)$, to lowest order 
the equation is $F+d a_1=0$. To lowest order, $[\mu,\mu]=0$
implies $dF=0$ so there is a solution $a_1$.
Then at each step one has to solve an equation of the
form $d(x)=b$ for given $b$ which is shown recursively to be
$d$-closed, as a consequence of $[\mu,\mu]=0$. 

As $\mu_1=\epsilon d+O(\epsilon^2)$, where $d$ is the Lie algebroid
differential, we have an $\mathbb R$\ndash linear map 
$p\colon H^0_{\pi,\epsilon}(N^*C)\to H^0_\pi(N^*C)$ sending $f_0+\epsilon f_1+\cdots
\in \mathrm{Ker}(\mu_1)|_{C^\infty(C)}$ to $f_0\in \mathrm{Ker}(d)$.
If $H^1_\pi(N^*C)$ vanishes there is a right inverse $\sigma\colon f_0\mapsto f$ 
to $p$ obtained
by solving recursively for $f_1,f_2,\dots$
the equation $d f=0$ with $f=f_0+\epsilon f_1+\cdots$.
By induction, at each step the equation is of the form 
$d f_j= C_j(f_0,\cdots,f_{j-1})$ with closed right-hand side and has a unique
solution $f_j\in K$. The ``quantisation map'' $\sigma$ is then extended by 
$\mathbb R[[\epsilon]]$\ndash linearity to an injective homomorphism 
$\sigma\colon H^0_{\pi}(N^*C)[[\epsilon]]\to H^0_{\pi,\epsilon}(N^*C)$ 
of $\mathbb R[[\epsilon]]$\ndash modules. 
The map $\sigma$ is surjective. Indeed,
if $f=f_0+\epsilon f_1+\cdots\in H^0_{\pi,\epsilon}(N^*C)$ then $f=\sigma(g_0+\epsilon g_1+\cdots)$ with $g_j$ recursively defined by $g_0=f_0$
$f-\sigma(g_0+\cdots+\epsilon^jg_j)=\epsilon^{j+1}g_{j+1}+
O(\epsilon^{j+2})$.

Since, by construction, $\sigma(f_0)=f_0+O(\epsilon)$, the
product $\mu_2$ induces a deformation of the
product on $H^0(N^*C)$ whose skew-symmetric part at 
order $\epsilon$ is, by Theorem \ref{t-Lowell}, is
the given Poisson bracket.
\end{prf} 

\medskip\noindent{\bf Remark.} If $H^2_\pi(N^*C)$ vanishes, one can also
construct
a quantisation of $C^\infty(\underline C)=H^0_\pi(N^*C)$ 
by the BRST method, see \cite{Bordemann}.

\medskip\noindent{\bf Remark.} From the explicit construction of the
$A_\infty$\ndash structure one sees that in some cases the obstruction
vanishes even if $H^2\neq 0$. For example if $\mathfrak h\subset\mathfrak g$
is an inclusion of finite dimensional real Lie algebras, 
then the subspace $C=\mathfrak h^\perp=(\mathfrak g/\mathfrak h)^*$
 of linear functions on $\mathfrak g$ vanishing on $\mathfrak
h$ is a coisotropic submanifold of the Poisson manifold
$\mathfrak g^*$ with Kostant--Kirillov bracket. In this case the
anomaly $\mu_0$ vanishes \cite{CaFe2003} even when 
$H^2_\pi(N^*C)=H^2_{\mathrm{Lie}}(\mathfrak h; C^\infty( \mathfrak h^\perp))
\neq 0$.

\section{The relative formality theorem}

\subsection{The Gerstenhaber algebra of multiderivations}
Let $A$ be a graded commutative algebra.
Recall that a derivation of degree $d$
of a graded algebra $A$ is a linear map $\varphi\colon
A\to A$ of degree $d$ such that $\varphi(ab)=\varphi(a)b
+(-1)^{d\cdot |a|}a\varphi(b)$, $a\in A^{|a|}$. 
Derivations form a graded left
$A$\ndash module $\mathrm{Der}(A)$ with a Lie bracket 
 $[\varphi,\psi]=\varphi\circ\psi-
(-1)^{|\varphi|\cdot|\psi|}\psi\circ\varphi$.
On the graded commutative algebra $S_A(\mathrm{Der}(A)[-1])$
(the graded symmetric algebra of the $A$-module $\mathrm{Der}(A)[-1]$)
 we then have a Gerstenhaber
structure, namely a (super) Lie bracket of degree $-1$ compatible
with the product. The Lie bracket is the extension of 
$[\ ,\ ]$  on
$\mathrm{Der}(A)$ to all of $S_A(\mathrm{Der}(A)[-1])$ 
 by the rule
\begin{equation}\label{e-Leibniz}
[\alpha\beta,\gamma]=\alpha[\beta,\gamma]
+(-1)^
{(\mathrm{deg}(\beta)+1)\cdot\mathrm{deg}(\gamma)}
[\alpha,\gamma]\beta.
\end{equation}
Here $\mathrm{deg}$ denotes the degree in the Lie
algebra of multiderivations
\[
\mathcal{T}(A)=S_A(\mathrm{Der}(A)[-1])[1],
\]
for which the Lie bracket has degree $0$ (and the
product degree $1$).
By definition $\mathrm{deg}(
\alpha_1\cdot\ldots\cdot\alpha_n)=\sum_{j=1}^n|\alpha_j|
+n-1$, for $\alpha_j\in\mathrm{Der}^{|\alpha_j|}(A)$.
The signs are then 
\[\alpha\beta=(-1)^{(\mathrm{deg}(\alpha)-1)(\mathrm{deg}(\beta)-1)}\beta\alpha,\qquad
[\alpha,\beta]=-(-1)^{\mathrm{deg}(\alpha)
\mathrm{deg}(\beta)}[\beta,\alpha].\]
 and the Jacobi identity is
\[
(-1)^{\mathrm{deg}(\alpha)\mathrm{deg}(\gamma)}[[\alpha,\beta],
\gamma]+\mathrm{cycl.}=0.
\]
To make contact with the Hochschild complex it will be useful to view multiderivations
as multilinear maps on $A$. First of all we have a map $\sigma\colon
S_A(\mathrm{Der}(A)[-1])\to
\oplus^j\wedge_A^j(\mathrm{Der}(A))[-j]$ given by
\[
\sigma(\varphi_1\cdot\ldots\cdot\varphi_j)
= (-1)^{\sum_{\alpha=1}^n(\alpha-1)|\varphi_\alpha|}
\varphi_1\wedge\dots\wedge\varphi_j.
\]
If $\mu\colon A^{\otimes n}\to A$ denotes the product in $A$, we have a map
$\tau\colon \wedge^j_A\mathrm{Der}(A)\to \mathrm{Hom}(A^{\otimes j},A)$:
\[
\tau(\varphi_1\wedge\cdots\wedge\varphi_j)=\mu\circ\mathrm{Alt}_j(\varphi_1\otimes\cdots
\otimes\varphi_j).
\]
The graded alternation map $\mathrm{Alt}_j\in\mathrm{End}(A^{\otimes j})$ is 
$(1/j!)\sum_{\sigma\in S_j}\mathrm{sign}(\sigma)\sigma$, for the natural action
(with Koszul signs)
of the symmetric group on the tensor algebra of the graded vector space $A$.

The {\em Hochschild\Ndash Kostant\Ndash Rosenberg map} is then the injective homomorphism
\begin{equation}\label{e-HKR}
\varphi_{\mathrm{HKR}}=\tau\circ \sigma\colon S_A(\mathrm{Der}(A)[-1])\mapsto 
\oplus_{j=0}^\infty\mathrm{Hom}_k(A^{\otimes j},A)
\end{equation}

\medskip

\noindent{\bf Remark.} In a more natural setting one would avoid the use of
exterior algebras and only use the symmetric algebras, which are always
graded commutative. The advantage would be that many signs would be simpler: for
example we would not have the sign coming from the definition of $\sigma$ and would
define the HKR map as a map to $\oplus\mathrm{Hom}(A[1]^{\otimes j},A[1])$. However
other things would become more exotic: for example, in this setting an associative
product (solution of the Maurer--Cartan equation) would be a map of degree 1 obeying
a ``graded associativity'' property, with signs.

\subsection{Fourier transform}\label{s-Fourier}
For the relative formality theorem two algebras and the corresponding Gerstenhaber
algebras are relevant. Let $C\subset M$ be a submanifold of a manifold $M$ which
we may assume to be the total space of a vector bundle $p:E\to C$ of rank $r$ with
$C$ embedded as the zero section.
Let
$A=\oplus_0^r A^j$, with $A^j=\Gamma(C,\wedge^j E)$ the graded
commutative algebra of sections of the exterior algebra of $E$. Let $B=\oplus B^j$, with
$B^0=\Gamma(C,S(E^*))$ and $B^j=0$ for $j\neq0$,  be the algebra of polynomial functions
on $E$, considered as a graded algebra concentrated in degree zero. Then the Gerstenhaber algebra $S_B(\mathrm{Der}(B)[-1])$
may be identified with the algebra of multivector fields on $M=E$ which are polynomial
along the fibres.

\begin{thm}\label{t-Fields}
The two Gerstenhaber algebras $S_A(\mathrm{Der}(A)[-1])$ and
$S_B(\mathrm{Der}(B)[-1])$ are canonically isomorphic up to choice
of sign.
In particular, the Lie algebras of multiderivations
$\mathcal{T}(A)$ and $\mathcal{T}(B)$ are canonically
isomorphic up to choice of sign.
\end{thm}  

Let us first suppose that $E$ is a trivial bundle $C\times V$ over
an open subset $C$ of $\mathbb{R}^n$ with coordinates
$x^1,\dots,x^n$. Let $\theta_1,\dots,
\theta_r$ be a basis of $V$. Then $A$ is freely generated by its
degree zero component
$A^0=C^\infty(C)$ and $\theta_\mu$ of degree 1 ($\mu=1,\dots,r)$. 
Then $\mathrm{Der}(A)$
is a free $A$\ndash module generated by 
$\xi_i=\partial/\partial x^i$ $(i=1,\dots,n$) of degree 0
and $\psi^\mu=\partial/\partial \theta_\mu$ $(\mu=1,\dots,r$)
of degree $-1$. Thus, as a graded 
algebra, $S_A(\mathrm{Der}(A)[-1])$
is the free graded commutative $A^0$\ndash algebra $A^0[\theta_\mu,
\psi^\mu,\xi_i]$ with generators $\theta_\mu$ of degree 1,
$\psi^\mu$ of degree $0$, and $\xi_i$ of degree $1$. 
The Lie bracket is defined by the relations
\begin{eqnarray*}
 [\xi_i,f]&=&\frac{\partial f}{\partial x^i},
\qquad [\psi^\mu,f]=[\theta_\mu,f]=0,  \qquad f\in C^\infty(C),\\
{} [\psi^\mu,\theta_\nu]&=&\delta^\mu_\nu,
\end{eqnarray*}
and the remaining brackets between generators vanish. 
Similarly, $B$ is generated by $C^\infty(C)$ and
the dual basis elements $y^\mu$ of $V^*$. We then
have $S_B(\mathrm{Der}(B)[-1])=A^0[y^\mu,\eta_\mu,\xi_i]$
with $\eta_\mu=\partial/\partial y^\mu$ of degree $1$.
The Lie bracket is 
\begin{eqnarray*}
[\xi_i,f]&=&\frac{\partial f}{\partial x^i},
\qquad [\eta_\mu,f]=[y^\mu,f]=0,  \qquad f\in C^\infty(C),\\
{}[\eta_\mu,y^\nu]&=&\delta_{\mu,\nu}.
\end{eqnarray*}
The isomorphism is then the isomorphism of graded commutative
algebras over $A^0$ that on generators is defined by
\[
\xi\mapsto\xi, \qquad \theta_\mu\mapsto -\eta_\mu,
\qquad \psi^\mu\mapsto y^\mu.
\]
To prove the theorem for general vector bundles, we
show that both algebras $S_A(\mathrm{Der}(A)[-1])$,
$S_B(\mathrm{Der}(B)[-1])$ are (non-canonically)
isomorphic, as graded commutative algebras, to
$R=\oplus_jR^j$, where
\[
R^j=\oplus_{p,q}\Gamma(C,\wedge^p E\otimes S^qE^*\otimes\wedge^{j-p}TC).
\]
The isomorphisms depend on the
choice of a connection $\nabla$ on $E$. This connection
induces a connection, also denoted by $\nabla$ on
$\wedge E$ and a dual 
connection $\nabla^*$ on $E^*$ and on $SE^*$.
Note that we have canonical inclusions of $A$ and
$B$ into $R$.
The isomorphism $j^\nabla_A\colon R\to S_A(\mathrm{Der}(A)[-1])$ sends
$A$ to $A$, $\psi\in\Gamma(C,S^1E^*)$ to the inner
multiplication $\iota_\psi\in\mathrm{Der}(A)$ and
$\xi\in\Gamma(C,TC)$ to $\nabla_\xi\in\mathrm{Der}(A)$.
The isomorphism $j^\nabla_B\colon R\to S_B(\mathrm{Der}(B)[-1])$ sends
$B$ to $B$, $\eta\in\Gamma(C,\wedge^1E)$ to $-\iota_\eta$
and $\xi\in\Gamma(C,TC)$ to $\nabla^*_\xi\in\mathrm{Der}(B)$.

\begin{lemma}
The composition of isomorphisms $j^\nabla_B\circ (j^\nabla_A)^{-1}$
is independent of the choice of connection $\nabla$ and respects the Lie brackets.
\end{lemma} 

\begin{prf}
It is sufficient to show that $j^\nabla_B\circ (j^\nabla_A)^{-1}$ sends
$\nabla'_\xi\in \mathrm{Der}(A)$ to $\nabla'{}^*_\xi$ for
any other connection $\nabla'$. The difference between
two connections is a 1-form with values in $\mathrm{End}(E)=
E\otimes E^*$. So we can write 
$\nabla'_\xi s=\nabla_\xi s +\sum a_i\wedge \iota_{b_i}s$
for some $b_i\in\Gamma(C,E^*)$, $a_i\in \Gamma(C,E)$ (depending
on $\xi$) and any $s\in A$. Thus the isomorphism maps 
$\nabla'_\xi$ to $\nabla^*_\xi-\sum  b_i\iota_{a_i}$ which
is precisely  $\nabla'{}^*_\xi$ on $\Gamma(C,E^*)$ and thus also on $B$.

With this result it is now easy to show that the isomorphism
respects the Lie bracket: as it is sufficient to prove this
locally, we may choose a connection which is locally the trivial
connection on the trivial bundle and use the local 
calculation above.
\end{prf}

\subsection{The Lie algebra of multidifferential operators}
Let $A$ be a graded commutative algebra, $C(A,A)$ the Hoch\-schild cochain
complex of $A$ (see \ref{s-Hoch}). The shifted complex $C(A,A)[1]$ is
a differential graded Lie algebra with respect to the Gerstenhaber bracket. It has a
subalgebra $\mathcal{D}(A)$ consisting of multidifferential operators, namely
sums of cochains of the form $(a_1,\dots,a_p)\mapsto \prod \varphi_i(a_i)$,
where $\varphi_i$ are compositions of derivations.
The HKR map $\mathcal{T}(A)\to \mathcal{D}(A)$
induces a homomorphism
of Gerstenhaber algebras on the
cohomology 
($\mathcal{T}(A)$ is considered as a complex
with zero differential).

\begin{lemma}
 If $A=\oplus_{j}\Gamma(C,\wedge^j E)$ for a vector bundle $E\to C$,
then the HKR map \eqref{e-HKR}, viewed as a map from $\mathcal{T}(A)$
with zero differential to $\mathcal{D}(A)$ with Hochschild differential,
induces an isomorphism on cohomology.
\end{lemma}

We prove this Lemma in the Appendix, see Lemma \ref{l-HKR}

\subsection{Completions} For our application, in the above construction we
take $E$ to be the normal bundle $NC$ to a submanifold $C\subset M$ with
vanishing ideal $I(C)=\{f\in C^\infty(M)\,|\, f|_C=0\}$. The relevant Lie
algebra is then $\mathcal{T}(M,C)=\varprojlim\mathcal{T}(M)/I(C)^n\mathcal{T}(M)$ 
of multivector fields on a formal neighbourhood of $C$. Let us fix an identification
of $NC$ with a tubular neighbourhood of $C$ as in \ref{s-PINF}. 
Introduce the graded commutative algebras $A=\oplus A^j$, $A^j=\Gamma(C,\wedge^j E)$
and $B=\Gamma(C,S(E^*))$ concentrated in degree $0$. Then $\mathcal{T}(M,C)$
may be viewed as the completion $\hat{\mathcal{T}}(B)=\varprojlim \mathcal{T}(B)/
I_B^n\mathcal{T}(B)$ of the $B$\ndash module $\mathcal{T}(B)$. Here $I_B$ is the
ideal $\Gamma(C,\oplus_{j>0}S^j(E^*))$ of $B$. 
Then there is a completion $\hat{\mathcal{T}}(A)$ defined by requiring the
isomorphism of  Theorem \ref{t-Fields} to extend to an isomorphism of the
completed Lie algebras. This completion is defined using the same ideal
$\Gamma(C,\oplus_{j>0}S^j(E^*))$ of $B$, which is now realised as the space of
$C^\infty(C)$\ndash multilinear multiderivations of $A$ with values in $C^\infty(C)$.
In both cases we have a Gerstenhaber algebra $G=\oplus_{j=0}^\infty G^j$
with nonnegative grading and a sequence of ideals (for the algebra structure)
$I^n\subset G^0$ such that
$[I^{n},G]\subset I^{n-1}$. It follows that the inverse limit 
is still a Gerstenhaber
algebra.

\begin{proposition}\label{p-Dalhousie}
 The image of a Poisson bracket on a formal neighbourhood of $C$ in
 $M$ under the (completed) isomorphism of Theorem \ref{t-Fields}
\[
\mathcal T^1(M,C)=\hat {\mathcal T}^1(B)\to\hat{\mathcal T}^1(A)
\]
is the $P_\infty$-structure of Proposition \ref{p-Pzeroinfty}
\end{proposition}

\begin{prf}
This can be proved in local coordinates using the trivial connection on trivial bundles
to describe the isomorphism: the result is that the components of the
$P_\infty$-structure are the Taylor expansion coefficients in the transverse coordinates
of the components of the Poisson bivector field.
\end{prf}

Now we need to find a completion of the Lie algebra of multidifferential operators
in such a way that the HKR map remains an isomorphism.
\dontprint{ Let $C^{p,q}(A,A)=
\mathrm{Hom}^p(A^{\otimes q},A)$, cf.~\ref{s-Hoch},  $C^{p,q}_{\mathrm{Diff}}(A,A)$
be the subspace of multidifferential operators, and $C^{m}_{\mathrm{Diff}}(A,A)=
\oplus_{p+q-1=m}C^{p,q}_{\mathrm{Diff}}(A,A)$. This space vanishes if $m<-1$ and if
$m=-1$, $C^{p,-p}_{\mathrm{Diff}}(A,A)$ has the local coordinate description
\[
(a_1,\dots,a_p)\mapsto f^{i_1\dots i_p}(x)
\frac{\partial{a_1}}{\partial\theta_{i_1}}\cdots\ 
\frac{\partial{a_p}}{\partial\theta_{i_p}}
\]
if $x_i$ are local commuting coordinates on $C$ and $\theta_i$ (of degree 1)
 are anticommuting generators of $A=\Gamma(C,\wedge E)$. Thus $\oplus_{p\geq 0}
C^{p,-p}_{\mathrm{Diff}}(A,A)$ is the image $\varphi_{\mathrm{HKR}}(B)$
 of $B$ by the HKR map. 
}

\medskip

\noindent{\bf Definition.} The completed Lie algebra of multidifferential
operators of $A=\Gamma(C,\wedge E)$ is
$\hat{\mathcal{D}}(A)=\oplus_{n} \hat{\mathcal{D}}^n(A)$,
where 
\[
\hat{\mathcal{D}}^n(A)=\prod_{p+q-1=n}\mathrm{Hom}^p(A^{\otimes q},A)
\] 
is the direct
product.

\medskip
The Gerstenhaber Lie bracket of two homogeneous elements $\phi=(\phi_{p,q})_{p+q-1=n}$,
$\psi=(\psi_{p,q})_{p+q-1=m}$
has $(p,q)$-component
\[
[\phi,\psi]_{p,q}=\sum [\phi_{p',q'},\psi_{p'',q''}],
\]
where the range of the sum is $p'+p''=p$, $q'+q''=q-1$, $p'+q'-1=n$, $p''+q''-1=m$
$q',q''\geq 0$, so we have a finite sum.
Clearly the HKR map extends naturally to an injective map $\hat{\mathcal{T}}(A)
\to\hat{\mathcal{D}}(A)$.

\begin{lemma}
The HKR map extends to a quasiisomorphism $\hat{\mathcal{T}}(A)\to\hat{\mathcal{D}}(A)$.
\end{lemma}

\begin{prf}
A cochain in $\hat{\mathcal{D}}(A)$ of degree $n$ is a sequence $\phi =(\phi_{p,q})_{p+q-1=n}$
with $\phi_{p,q}\in\mathrm{Hom}^p(A^{\otimes q},A)$. As the Hochschild differential
only shifts $q$, $\phi$ is a cocycle if and only if $b\phi_{p,q}=0$ for all $p,q$. By
the HKR theorem for $\mathcal{D}(A)$, $\phi_{p,q}=\psi_{p,q} \mod b \mathrm{Hom}^p(A^{\otimes{(q-1))}},A)$, for a unique $\psi_{p,q}$  in the image of the HKR map. 
This implies that $\phi=\psi+$exact for a unique $\psi$ in the image of the
extension of the HKR map to $\hat {\mathcal{T}}(A)$.
\end{prf}

\subsection{A graded version of Kontsevich's theorem}
Let $E\to C$ be a vector bundle on a smooth manifold $C$ and
$A=\oplus_{j=0}^\infty\Gamma(C,\wedge^j E)$. Let $\mathcal{T}(A)=
S_A(\mathrm{Der}(A)[-1])[1]$ be the differential graded Lie algebra of
multiderivations of $A$ with Schouten\Ndash Nijenhuis bracket. Let
$\mathcal{D}(A)$ be the subcomplex of the shifted Hochschild complex $C(A,A)[1]$
consisting of multidifferential operators, with Lie algebra structure
given by the Gerstenhaber bracket. 
In the language of supermanifolds, $A$ is (by definition) 
the algebra of smooth functions on the supermanifold $\Pi E^*$ obtained
from $E$ by changing the parity of the fibres and $\mathcal{T}(A)$ is
the Lie algebra of multivector fields on $\Pi E^*$.

\begin{thm}\label{t-Holmes}  There exists
an $L_\infty$\ndash quasiisomorphism $U\colon \mathcal{T}(A)\to \mathcal{D}(A)$
whose first order term $U_1$ is the Hoch\-schild\Ndash Kostant\Ndash Ro\-sen\-berg map
\eqref{e-HKR}
\[
\gamma_1\cdots\gamma_p\mapsto\frac1{p!}\mu\circ \sum_{\sigma\in S_p}
(-1)^{\sum_{i<j,\sigma(i)>\sigma(j)}
\mathrm{deg}(\gamma_i)\mathrm{deg}(\gamma_j)} \gamma_{\sigma(1)}\otimes\cdots\otimes\gamma_{\sigma(p)}
\]
\end{thm}

This theorem is implicitly stated in \cite{Kontsevich}. We give a construction of an $L_\infty$\ndash
quasiisomorphism in the Appendix. Composing this $L_\infty$\ndash isomorphism with the Lie algebra
isomorphism of \ref{s-Fourier}, we obtain the proof of Theorem \ref{t-rft}.

We need a completed version of this theorem. For this  the following property
of the $L_\infty$\ndash morphism is important. Suppose $\gamma_1,\dots,\gamma_n\in 
\mathcal{T}(A)$ are of order $p_1,\dots,p_n$, i.e., $\gamma_i\in S^{p_i}(\mathrm{Der}(A)[1])[-1]$. Let $f_1,\dots,f_m\in A$. Then $U_n(\gamma_1,\dots,\gamma_n)(f_1,\dots,f_m)$
vanishes unless
\[
\sum_{i=1}^np_i=2n+m-2.
\]
This condition expresses the fact that the degree of the differential forms appearing
in the definition of the weights entering $U_n$ coincides
with the dimension of the configuration spaces over which these differential forms
are integrated.
It follows that for given $n$ and $m$ there are only finitely many values
of $(p_1,\dots,p_n)$ giving a non-trivial contribution and thus all Taylor
components $U_n$ are well-defined on $\hat{\mathcal{T}}(A)$. We thus obtain:

\begin{thm}\label{t-Holmes2}  There exists
an $L_\infty$\ndash quasiisomorphism $U\colon \hat{\mathcal{T}}(A)\to \hat{\mathcal{D}}(A)$
whose first order term $U_1$ is the Hoch\-schild\Ndash Kostant\Ndash Ro\-sen\-berg map.
\end{thm}

\subsection{The relative formality theorem}

\begin{thm}\label{t-rft}
Let $C\subset M$ be a submanifold of a smooth manifold $M$ with vanishing ideal $I(C)$,
$A=\Gamma(C,\wedge NC)$ the graded commutative algebra of sections of
the exterior algebra of the normal bundle, $\mathcal{T}(M)=\Gamma(M,\wedge TM)$ 
the DGLA of multivector fields with Nijenhuis\Ndash Schouten
bracket and zero differential, $\mathcal{T}(M,C)=\varprojlim \mathcal{T}(M)/I(C)^n\mathcal{T}(M)$ 
the DGLA of multivector fields in an infinitesimal neighbourhood of $C$.
 Then there is an $L_\infty$\ndash quasiisomorphism
\[
U\colon\mathcal{T}(M,C)\to \hat{\mathcal{D}}(A),
\]
whose first order term $U_1$ is the composition
\[
\mathcal{T}(M,C)\simeq \hat{\mathcal{T}}(B)
\to \hat{\mathcal{T}}(A)
\stackrel{\varphi_\mathrm{HKR}}
{\longrightarrow}\hat{\mathcal{D}}(A).
\]
where the middle arrow is the Fourier transform isomorphism of Theorem \ref{t-Fields}
and $\varphi_\mathrm{HKR}$ is the HKR map \eqref{e-HKR}.
\end{thm}

Theorem \ref{t-rft} follows from Theorem \ref{t-Fields} and Theorem \ref{t-Holmes2}.

The $L_\infty$\ndash quasiisomorphism induces a bijection between deformation functors
(see \cite{Kontsevich}). In particular, a Poisson bivector field $\pi$
on $M$ defines a solution $\epsilon\pi$
of the Maurer Cartan equation $[\epsilon\pi,\epsilon\pi]=0$ 
in the pronilpotent Lie algebra $\mathcal{T}(M,C)\otimes\epsilon\mathbb{R}[[\epsilon]]$.
This solution is mapped by $U$ to a solution $\mu$ of the Maurer\Ndash Cartan equation
$2b\mu+[\mu,\mu]=0$ in $\epsilon\hat{\mathcal{D}}(A)[[\epsilon]]$, i.e., a
deformation of the product on $A$ as an $A_{\infty}$\ndash algebra. This proves
Theorem \ref{t-azeroinfty}.

\appendix
 
\section{Formality theorem for supermanifolds}
This section contains a graded version of
Kontsevich's formality theorem, stating that the
differential graded Lie algebra of 
multidifferential operators
on a graded super vector spaces is 
$L_\infty$\ndash quasiisomorphic to its cohomology, the graded Lie
algebra of multivector fields.
The proof is the same as Kontsevich's proof in
the case of ordinary vector spaces. Our contribution
is to write all signs and develop a formulation
in which the signs (which
are already non\ndash trivial in the case of
ordinary vector spaces) appear in a transparent
way. 

\subsection{Notations and conventions}\label{N&C} We work in the tensor category of
graded vector spaces over a field or more generally of
graded left modules over a graded commutative ring $R$ with unit.
All graded modules shall be meant to be $\mathbb Z$\ndash graded and shall
be considered as super vector spaces with the induced
$\mathbb Z/2\mathbb Z$-grading. The word super shall usually
be omitted. Thus
$R$ is a $\mathbb Z$\ndash graded commutative ring $R=\oplus_{j\in\mathbb Z}R^j$
and an object is a $\mathbb Z$\ndash graded left $R$\ndash module $V=\oplus V^j$.
Morphisms from $V$ to $W$ form a $\mathbb Z$\ndash graded left $R$\ndash module
$\mathrm{Hom}(V,W)=\oplus_d \mathrm{Hom}^d(V,W)$. We denote by $|a|$
the degree of a homogeneous element $a$. The Koszul sign rule
holds. Thus a homogeneous morphism $\phi\in \mathrm{Hom}(V,W)$ is an
additive map obeying $\phi(rv)=(-1)^{|r||\phi|}r\phi(v)$, $r\in R$,
$v\in V$; the tensor product $V\otimes W$ of objects is defined as the
quotient of the tensor product over $\mathbb Z$ by  the relation
$rv\otimes w=(-1)^{|r||v|}v\otimes rw$, $r\in R$; the tensor products
of morphisms $\phi\in\mathrm{Hom}(V,V')$,
$\psi\in\mathrm{Hom}(W,W')$ 
is $\phi\otimes\psi(v\otimes w)=(-1)^{|\psi||v|}\phi(v)\otimes
\psi(w)$, $v\in V$, $w\in W$. 

For a graded $R$\ndash module $V$ let $V[n]$ be the graded $R$\ndash module such that
$V[n]^j=V^{n+j}$. We have a tautological map (the identity) $s^n\colon V[n]\to V$
of degree $n$.

 We often denote by $(v_1,\dots,v_n)$ the element
$v_1\otimes\cdots\otimes v_n\in V^{\otimes n}=V\otimes\cdots\otimes V$.
The symmetric group $S_n$ acts on $V^{\otimes n}$ with signs: so
the transposition $s_i=(i,i+1)$ acts as 
\[s_i(v_1,\dots,v_n)=(-1)^{|v_i||v_{i+1}|}
(v_1,\dots,v_{i+1},v_i,\dots,v_n).\]
 The product of symmetric groups $S$ acts on the 
tensor algebra $T(V)=\oplus_{n\geq 0} V^{\otimes n}$ (with $V^{\otimes 0}=R$) and
we have the algebras of coinvariants for the ordinary action $S(V)=T(V)/(x-\sigma x),\sigma\in S$
and for the alternating action $\wedge(V)=T(V)/(x-\mathrm{sign}(\sigma)\sigma x,\sigma\in S$.

\subsection{Tensor coalgebras}\label{s-free}
Let $V$ be a free graded $R$\ndash module over a commutative unital
ring $R$ ($=\mathbb R$ or $\mathbb R[[\epsilon]]$
in our application).
We set $V^{\otimes 0}=R$ and $V^{\otimes j}=V\otimes\cdots\otimes V$. 
The  graded counital tensor coalgebra generated by $V$ is the graded $R$\ndash module
$F(V)=\oplus_{j\geq 0}V^{\otimes j}$ with coproduct
\[
\Delta(\gamma_1,\dots,\gamma_n)=\sum_{j=0}^n(\gamma_1,\dots,\gamma_j)\otimes(\gamma_{j+1},\dots,\gamma_n).
\]
In the first and last term we have $()=1\in R$. The counit is the
canonical projection onto $V^{\otimes 0}=R$.
The spaces of invariant tensors $I_n(V)=\{v\in V^{\otimes n}\,|\,\sigma v=v
\;\forall \sigma\in S_n\}$ form a commutative
sub\ndash coalgebra $C(V)=\oplus_{n\geq0} I_n(V)$
of $F(V)$. It is the symmetric coalgebra generated by $V$.
The quotient $F^0(V)$ by the coideal $R=V^{\otimes 0}$ can be
described as the coalgebra $\oplus_{j\geq1} V^{\otimes j}$ without counit
and whose coproduct is given by the formula above without the first
and last term. Similarly we have the coalgebra $C^0(V)=C(V)/R$. The coalgebras
$C^0(V)$, $F^0(V)$ are freely generated by $V$. For $C^0(V)$ this means that
if $C$ is a cocommutative coalgebra without counit so that, for each $x\in C$, 
the iterated coproduct $\Delta^n(x)$ vanishes for $n$ large enough, 
every linear map $U\colon C\to V$
is uniquely the composition of a map of coalgebras $\bar U\colon C\to C^0(V)$ with the
canonical projection $p_1:C^0(V)\to V$ on the first direct summand $V=I_1(V)$.
The formula for the composition of $\bar U$ with the canonical projection $p_n:C^0(V)\to I_n(V)$
is 
\[
p_n\circ\bar U(x)=U\otimes\cdots\otimes U(\Delta^{n}(x)). 
 \]
For $F(V)$, $C(V)$ we will need the following infinitesimal version of this fact in a special case.

\begin{lemma}\label{l-Alighieri}
Let $Q\colon F(V)\to V$ be a linear map, $p_n:F(V)\to V^{\otimes n}$ the canonical projection onto the
$n$th summand. Then there is a unique coderivation $\bar Q\colon F(V)\to F(V)$ such that $p_1\circ\bar Q =Q$.
The same holds for in the cocommutative case for $C(V)$ with projections $p_n\colon C(V)\to I_n(V)$.
\end{lemma}

The formula for the components of 
$\bar Q_n=\bar Q|_{V^{\otimes n}}$ (or $\bar Q|_{I_n(V)}$) in terms of the components of $Q$ is
\[
\bar Q_n=\sum_{m=0}^n\sum_{l=0}^{n-m}1^{\otimes l}\otimes Q_m\otimes 1^{\otimes {n-m-l}},
\]
where $1^{\otimes l}$ denotes the identity on $V^{\otimes l}$. 

Closely related to $C(V)$ is the shuffle coalgebra $S(V)$,
the $R$\ndash module of coinvariants $\oplus V^{\otimes j}/\{\sigma v-v,
\sigma\in S_j\}$ with shuffle coproduct
\[
\Delta_{\mathrm{sh}}(\gamma_1,\dots,\gamma_n)=\sum_{p+q=n}
\sum_{(p,q)-\mathrm{shuffles}}\pm
(\gamma_{\sigma(1)},\dots,\gamma_{\sigma(p)})
\otimes
(\gamma_{\sigma(p+1)},\dots,\gamma_{\sigma(n)}).
\]
The sum is over permutations such that $\sigma(1)<\cdots<\sigma(p)$,
$\sigma(p+1)<\cdots<\sigma(n)$ with sign
\begin{equation}\label{e-sign}
\pm=
\eps(\sigma,\gamma_1,\dots,\gamma_n)=(-1)^{\sum_{i<j,\sigma(i)>\sigma(j)}
|\gamma_i|\cdot|\gamma_j|}.
\end{equation}
The map $S(V)\to C(V)$ sending $(\gamma_1,\dots,\gamma_n)$ to 
$\sum_{\sigma\in S_n}\pm(\gamma_{\sigma(1)},\dots,\gamma_{\sigma(n)})$
is then an isomorphism of coalgebras.

Although the language of coalgebras is technically convenient, it
is better for the intuition to think in terms of the dual algebras.
If $R$ is the field of real or complex numbers,
the dual space $C(V)^*=\mathrm{Hom}(C(V),R)$ is the algebra of
jets at zero of functions on the (super)manifold $V$ and $C^0(V)^*$ the
subalgebra of functions vanishing at $0$ and a coderivation of $C(V)$ is
a formal vector field. The algebras $F(V)^*$, $F^0(V)^*$
are the corresponding non\ndash commutative analogues.

\subsection{The Hochschild complex of a graded algebra}\label{s-Hoch}
Let $A=\oplus_{j\in\mathbb Z} A^j$ be a graded associative
algebra with unit over a field $k$.  
The Hochschild complex $C(A,A)$ 
with values in $A$, is the
complex
$C(A,A)=\oplus_{n} C^{n}(A,A)$
where
\[
C^{n}(A,A)=\oplus_{m+d=n}C^{d,m}(A,A),
\qquad C^{d,m}(A,A)=\mathrm{Hom}^d(A^{\otimes m},A) 
\]
The Hochschild differential of $\phi\in C^{|\phi|,m}(A,A)$ is
\begin{eqnarray}\label{e-Hoch}
b\phi(a_1,\dots,a_{m+1})&=&(-1)^{|\phi||a_1|}
a_1\phi(a_2,\dots,a_{m+1})\notag\\ 
&&+\sum_{j=1}^m(-1)^j
\phi(a_1,\dots,a_ja_{j+1},\dots,a_{m+1})\\
&&+(-1)^{m+1}
\phi(a_1,\dots,a_m)a_{m+1}.\notag
\end{eqnarray}
The shifted Hochschild complex $C(A,A)[1]$ is a differential graded Lie 
algebra whose Lie bracket is (a graded version of) the 
Gerstenhaber bracket: let $\phi\in C^{|\phi|,m_1}(A,A)$,
$\psi\in C^{|\psi|,m_2}(A,A)$.
\begin{equation}\label{e-GB}
[\phi,\psi]_\mathrm{G}=\phi\bullet\psi-
(-1)^{(|\phi|+m_1-1)(|\psi|+m_2-1)}\psi\bullet\phi,
\end{equation}
with Gerstenhaber product\footnote{In the ungraded case, 
this product differs by a factor
$(-1)^{(m_1-1)(m_2-1)}$ from the product defined in \cite{Gerstenhaber}. With
our convention we obtain a more standard bracket on multivector fields}
\[
\phi\bullet\psi=(-1)^{(|\psi|+m_2-1)(m_1-1)}\sum_{l=0}^{m_1-1}
(-1)^{l(m_2-1)}\phi\circ (1^{\otimes l}\otimes
\psi\otimes 1^{\otimes (m_1-1-l)}).
\]
We also have the cup product on $C(A,A)$: $\phi_1\cup\phi_2
=\mu\circ \phi_1\otimes \phi_2$, where $\mu$ is the product
in $A$. 

The simplest way to prove that the Gerstenhaber bracket
is a Lie bracket is to use its interpretation as 
a commutator of coderivations \cite{Stasheff93}. This also explains 
the origin of the signs:
a coderivation of degree $d$
of a coalgebra $C$ is a linear endomorphism
$\phi\in \mathrm{Hom}^d(C,C)$ 
obeying $\Delta\circ\phi=
(\phi\otimes 1+1\otimes \phi)\circ\Delta$. 
Coderivations of $C$ form a graded Lie algebra
$\mathrm{Coder}(C)$ with respect to the graded commutator
$\phi\circ\psi-(-1)^{|\phi||\psi|}\psi\circ\phi$.
Let $C=F(A[1])$. Then
a coderivation $\phi$ of $C$ 
is uniquely determined by its composition $p_1\circ\phi$ with the canonical projection onto $A[1]$
and every map $F(A[1])\to A[1]$ extends to a coderivation.
Thus we can identify derivations of $C$ 
with maps $F(A[1])\to A[1]$. Under this identification,
the Lie bracket is
\[
[\phi,\psi]=\sum_{l=0}^{m_1-1}\phi\circ (1^{\otimes l}\otimes
\psi\otimes 1^{\otimes (m_1-1-l)})-(-1)^{|\phi||\psi|}(
\phi\leftrightarrow \psi).
\]
 if $\phi\in\mathrm{Hom}(A[1]^{\otimes m_1},A[1])$ and
$\psi\in\mathrm{Hom}(A[1]^{\otimes m_2},A[1])$. Let
$s\colon A[1]\to A$ be the tautological map (of degree $1$)
and introduce $\tilde \phi$ by
 $\phi=s^{-1}\circ
\tilde\phi \circ(s\otimes \cdots\otimes s).$ Then the
Gerstenhaber bracket is
\[
[\tilde\phi,\tilde\psi]_G=\widetilde{[\phi,\psi]},
\] 
and the signs are obtained by the Koszul rule when
letting the maps $s$ go past $\phi,\psi$ and other maps $s$.

The Hochschild differential can also be expressed in terms of
the bracket: let $\mu\colon A\otimes  A\to A$ denote the
product in $A$. Then the associativity is the relation $[\mu,\mu]_G=0$. It 
follows that $[\mu,\cdot] $ is a differential and indeed
\[
b\phi=(-1)^{|\phi|}[\mu,\phi]=-[\phi,\mu].
\]
The cohomology of $C(A,A)$ is denoted by $\mathit{HH}(A,A)$. The Gerstenhaber
bracket induces a graded Lie algebra structure on $\mathit{HH}(A,A)[1]$.
 In terms of homological algebra, $\mathit{HH}(A,A)=\mathrm{Ext}_{A-A}(A,A)$ is the $\mathrm{Ext}$
group of $A$ in the category of $A-A$\ndash bimodules over the graded algebra $A$. Indeed
$C(A,A)\simeq\mathrm{Hom}_{A-A}(B(A),A)$, where $B(A)=\oplus_j (A\otimes A^{\otimes j}\otimes A)$
is the bar resolution, with degree assignment
\[|a\otimes a_1\otimes\cdots\otimes a_j\otimes b|=|a|+\sum_{i=1}^j|a_i|+|b|-j.\]
The sign in \eqref{e-Hoch}
comes from the Koszul rule for morphisms $\phi\colon M\to N$ of $A-A$\ndash bimodules:
\[
\phi(a\,m\,b)=(-1)^{|a||\phi|}a\,\phi(m)\,b, \qquad a,b\in A, m\in M
\]
\subsection{HKR cocycles}
Let $A$ be a graded commutative algebra over a field $k$ of characteristic zero, 
$C(A,A)=\oplus_{n\geq0}
\mathrm{Hom}(A^{\otimes n},A)$ the Hoch\-schild cochain complex of $A$. 

A cochain $\phi\in C(A,A)$ is called an {\em HKR-cocycle} if (i) the map
$a\mapsto\phi(a_1,\dots,a_{n-1},a)$ is a derivation 
of $A$ for any homogeneous $a_1,\dots,a_{n-1}\in A$
and (ii) $\phi$ is alternating in the graded sense, i.e., 
\[
\phi(a_1,\dots,a_i,a_{i+1},\dots,a_n)
=-(-1)^{|a_i||a_{i+1}|}
\phi(a_1,\dots,a_{i+1},a_i,\dots,a_n).
\]
It is easy to check that cochains  obeying (i) and (ii)
are indeed cocycles. They thus form a subcomplex $C_{\mathrm{HKR}}(A,A)$ with
zero differential that can be identified with $\wedge_A\mathrm{Der}(A)$
via the map $\phi_1\wedge\cdots\wedge\phi_n\mapsto \mu\circ\mathrm{Alt}(\phi_1\otimes
\cdots\otimes\phi_n)$ for derivations $\phi_j$; here $\mu\colon A^{\otimes n}\to A$ is
the product. 
 
A cochain $\phi\in C(A,A)$ is called {\em multidifferential operator} if
it is a sum of terms of the form $a_1\otimes\cdots\otimes a_n\mapsto
D_1(a_1)\cdots D_n(a_n)$, for some differential operators (compositions of
derivations) $D_i$. Multidifferential operators form a subcomplex $C_\mathrm{Diff}(A,A)$
of $C(A,A)$ containing $C_{\mathrm{HKR}}(A,A)$. 

Let $N$ be a graded supermanifold. For us this means that the ground field $k$ is
$\mathbb R$ (or $\mathbb C$) and $N=E^*$ is
the total space of the dual of a graded vector bundle $E\to C$ so that 
$C^\infty(N)=\Gamma(C,S(E))$, where $S(E)\to C$ 
is the graded symmetric algebra of $E$
(thus $S(E)=S(E^{\mathrm{even}})\otimes\wedge(E^{\mathrm{odd}})$).

\begin{lemma}\label{l-HKR} If $A=C^\infty(N)$ for a graded supermanifold $N$
the HKR map $\wedge_A\mathrm{Der}(A)\simeq C_\mathrm{HKR}(A,A)\hookrightarrow
C_\mathrm{Diff}(A,A)$ is a quasiisomorphism of complexes.
\end{lemma}

In the ungraded case $(E=0)$ a version of this theorem can be found in
\cite{Vey}, see also \cite{Kontsevich}. It is an analogue for smooth functions of the
original HKR theorem \cite{HKR}, which deals with regular affine algebras.

The proof is the same as the proof in the ungraded case (see \cite{Kontsevich},
4.4.1.1) but with some twists. First one uses the filtration by the
total order of multidifferential operators to pass to the associated
graded complexes of principal symbols. These complexes are sections
of vector bundles and the differential is 
$C^\infty(N)$\ndash linear, so the problem is reduced to proving a
version of the HKR theorem for each fibre. If $T=T_xM\oplus E_x$ 
is a tangent space to $N$, the complex of
principal symbols at a point $x\in M$
is $\oplus_{n\geq 0}S(T)^{\otimes n}$ with degree
assignment $|D_1\otimes\cdots\otimes D_n|=\sum |D_i|+n$, $D_i\in S(T)$.
An element of $S(T)$, considered as a differential operator with
constant coefficient, defines a linear function on the algebra
$S(T^*)$ of polynomial functions on the graded vector space $T$. 
Thus we obtain an embedding
$\oplus_{n\geq 0}S(T)^{\otimes n}\to \mathrm{Hom}_k (S(T^*)^{\otimes
  n},k)$ as a subcomplex. The differential on
$\mathrm{Hom}_k(S(T^*)^{\otimes n},k)$ is
\begin{eqnarray*}
  d\varphi(f_1,\dots,f_{n+1})&=&(-1)^{|f_1||\varphi|}
  \epsilon(f_1)\varphi(f_2,\dots,f_{n+1})\\
  &&+
  \sum_{j=1}^n(-1)^{j}\varphi(f_1,\dots,f_jf_{j+1},\dots,f_{n+1})
  \\ &&-(-1)^{n}\varphi(f_1,\dots,f_n).
\end{eqnarray*}
Here $\epsilon(f)=f(0)$.

\begin{lemma} The map of complexes 
  $(\wedge T,0)\to (\oplus_{n\geq 0} S(T)^{\otimes n},d)$ sending
  $t_1\wedge\cdots\wedge t_n$ to $\mathrm{Alt}(t_1\otimes\cdots\otimes
  t_n)$ is a quasiisomorphism.
\end{lemma}

\noindent{\em Proof:} Let $S=S(T^*)$ and view both complexes as subcomplexes
of $(\oplus_j C^j(S,k),d)$, 
\[
C^j(S,k)=\prod_{p+q=j}\mathrm{Hom}^p(S^{\otimes q},k).
\]
 In particular $(\wedge T,0)$ is
identified with the subcomplex $C_{\mathrm{HKR}}(S,k)$ consisting of 
cochains obeying (i) and (ii) above. We show that the embedding of
this subcomplex is a quasiisomorphism.  Since the complex $C(S,k)$ is a direct
product of subcomplexes consisting of multidifferential operators of 
fixed total order, it then follows that the same statement holds if we
replace $C(S,k)$ by its subcomplex $S(T)$.

To compute the cohomology of $C(S,k)$ we first notice that it is
 $\mathrm{Ext}_{\mathrm{mod}-S}(k,k)$, where $k$ is considered as a right $S$-module via
$\epsilon$ and has a free resolution $\cdots\to S^{\otimes 3}\to
S^{\otimes 2}\to S\to k$ with differential 
\[
a_1\otimes\cdots\otimes a_{n+1}\to \epsilon(a_1)a_2\otimes\cdots\otimes a_{n+1}-
\sum_{i=1}^{n} (-1)^i a_1\otimes\cdots\otimes a_ia_{i+1}\otimes \cdots\otimes a_{n+1},
\]
inducing the differential above. This
Ext group (in the category of graded right $S$\ndash modules)
can be computed using a graded version of the Koszul
resolution of the $S$-module $k$: let $v_1,\dots,v_n$ be a
homogeneous basis of the graded vector space $T^*$ so that $S$ is the
graded polynomial algebra $k[v_1,\dots,v_n]$.  Let
$K=S[u_1,\dots,u_n]=k[v_1,\dots,v_n,u_1,\dots,u_n]$ be the
differential graded commutative algebra with $u_i$ of degree $|v_i|+1$
and differential $\partial$ such that $\partial u_i=v_i, \partial
v_i=0$. Then $K$ is a free $S$-module and the map $(K,\partial)\to
(k[u_1,\dots,u_n],0)$ is a quasiisomorphism. The proof of the latter
statement is similar to the one in the ungraded case (see, e.g.,
\cite{MacLane} VII.2): since $K$ is the (graded) tensor product of
algebras $k[v_i,u_i]$, it is sufficient to check this for $n=1$. In
this case there is a homotopy $h:K\to K$ of degree $1$ obeying
$h\partial+\partial h=\mathrm{id}-\epsilon$ from which the claim
follows immediately: if $x=v_i$ is even, $h(x^p)=x^{p-1}u$ $(p\geq1)$,
$h(1)=0=h(x^pu)$; if $x$ is odd, $h(u^p)=u^{p+1}/(p+1)$, $h(xu^p)=0$.
Thus $\mathrm{Ext}_S(k,k)$ is the cohomology of
\[
\mathrm{Hom}_S(K,k)=\mathrm{Hom}_k(k[u_1,\dots,u_n],k).
\]
Since the induced differential vanishes identically ($v_i$ acts by zero on $k$),
we obtain
\[
\mathrm{Ext}^j_{S(T^*)}(k,k)=\mathrm{Hom}^j_k(k[u_1,\dots,u_n],k)
\]
This space may be identified with $\wedge T$. To find the map, we need to
write the map between the two resolutions, which is known to exist from 
abstract nonsense. Its explicit expression is
\[
u_{i_1}\cdots u_{i_n}a\mapsto (-1)^{\sum_{\alpha=1}^n(\alpha-1)d(i_\alpha)}
\mathrm{Alt}(v_{i_1}\otimes\cdots\otimes v_{i_n})\otimes a,\qquad a\in S,
\]
where $d(i)=|v_i|$ is the degree of $v_i$. The claim of the Lemma then follows
from the fact that the restriction to 
$C^j_{\mathrm{HKR}}(S,k)$
of the dual map $\oplus_{p+q=j}
\mathrm{Hom}^{ p}(S^{\otimes q},k)\to
\mathrm{Hom}^j(k[u_1,\dots,u_n],k)$  is an isomorphism.

\subsection{Expressions in local coordinates}

Let $A=S(V)$ be the algebra of polynomial functions on a finite dimensional 
graded real vector space $V^*$. 
If $(x_i)_{i=1}^d$ is a homogeneous basis of degrees $\epsilon_i=|x_i|$, 
then $S(V)$ is
the free graded commutative algebra $\mathbb R[x_1,\dots,x_d]$ generated by the $x_i$'s.
The Gerstenhaber algebra $S_A(\mathrm{Der}(A)[-1])$ may then be identified
$\tilde A=S(V\oplus V[1]^*)=k[x_1,\dots,x_n,\theta_1,\dots,\theta_n]$ where
$\theta_i$ is the dual basis with degrees $|\theta_i|=1-\epsilon_i$. 
 Write a  general element of
$\tilde A^{|\gamma|}=\tilde A[1]^{|\gamma|-1}$ with the summation convention as
\[
\gamma=
\gamma^{i_1\cdots i_m}
\theta_{i_1}\cdots\theta_{i_m},\quad |\gamma|=|\gamma^{i_1\cdots i_m}|+m-\sum\eps_{i_\alpha}.
\]
with $\gamma^{\dots,i,j,\dots}=(-1)^{(1-\eps_i)
(1-\eps_j)}\gamma^{\dots,j,i,\dots}\in A$.
The  HKR map  is then 
\[
\gamma\mapsto
(-1)^{\sum_{\alpha=1}^m(\alpha-1)\eps_{i_\alpha}}
\gamma^{i_1\cdots i_m}
\mu\circ(\partial_{i_1}\otimes\cdots\otimes
\partial_{i_m}).
\]
Here $\mu(f_1,\dots,f_m)=f_1\cdots f_m$ is the product
in $A$.

The Lie algebra structure on $S(V)\otimes S(V[1]^*)=S(V\oplus V[1]^*)$ induced by the HKR homomorphism may be
understood geometrically as the Poisson structure on the functions on the
degree\ndash shifted cotangent bundle $T^*[1]M$ of the supermanifold $M=V^*$ 
with its canonical odd symplectic structure. 
Here is the explicit description.
On $\tilde A$ there is a Poisson bracket of degree $-1$:
\[
[\gamma_1,\gamma_2]=\sum_{i=1}^d
\left(\gamma_1\!
\stackrel{\leftarrow}\partial_{\theta_i}
\stackrel{\rightarrow}\partial_{x_i}\!\!
\gamma_2
-
\gamma_1\!
\stackrel{\leftarrow}\partial_{x_i}
\stackrel{\rightarrow}\partial_{\theta_i}\!\!
\gamma_2
\right).\]
The operator $\stackrel{\leftarrow}\partial_x$ is
the right partial derivative with respect to $x$ acting
on the argument on its left as a right derivation
$(ab)\!\!\stackrel{\leftarrow}{\partial_x}=a(b\!\!\stackrel{\leftarrow}\partial_x\!)+(-1)^{|b||x|}(a\!\stackrel{\leftarrow}\partial_x\!)b$. The left derivatives $\stackrel{\rightarrow}\partial_x$ are defined in the same way and act to the right.
In more standard notation (using only left derivatives),
\[
[\gamma_1,\gamma_2]=\sum_{i=1}^d
(-1)^{(1-\eps_i)(|\gamma_1|-1)}
\partial_{\theta_i}\gamma_1
\partial_{x_i}
\gamma_2
-
(-1)^{\eps_i(|\gamma_1|-1)}
\partial_{x_i}\gamma_1
\partial_{\theta_i}\gamma_2.
\]
Shifting degrees we obtain thus a Lie algebra
$\mathcal{T}(A)=\tilde A[1]$, the graded Lie algebra of {\em multivector fields}.

\dontprint{\medskip

\noindent{\bf Definition.} A {\em multidifferential
operator}  is a cochain in $C(A,A)$ which
can be written as a sum of terms of the form
\[
 (a_1,\dots,a_m)\mapsto \phi_1(a_1)\cdots\phi_n(a_n),
\]
where $\phi_i\in A[\partial_{x_1},\dots,\partial_{x_d}]$
are differential operators. 

\medskip

Multidifferential operators form a differential graded Lie 
subalgebra $\mathcal{D}(A)$ of $C(A,A)$. The HKR homomorphism maps is then a quasiisomorphism
$\mathcal{T}(A)\to \mathcal{D}(A)$. The advantage of replacing $C(A,A)$ by
the subalgebra $\mathcal{D}(A)$  is that the construction also works if we replace
$A$ by the algebra of formal power series in $x_i$.
}

\subsection{$Q$\ndash manifolds and $L_\infty$\ndash algebras}
We use the language of (formal, pointed) $Q$\ndash manifolds. 
Let $V=\oplus_{j\in\mathbb Z}V^j$ be a graded real vector space.
Let $C^0(V)=\oplus_{j=1^\infty} I_j(V)$
be the free cocommutative coalgebra without counit generated by $V$. 
Its dual is the algebra of functions in an infinitesimal
neighbourhood of $0$ in $V$. 
A (formal, pointed) $Q$\ndash manifold is a graded vector
space $V$ with a coderivation $Q$ of $C^0(V)$ of degree $1$ 
obeying $[Q,Q]=0$. Dually, $Q$ may be thought of as a
vector field of degree 1 defined on a formal neighbourhood
of $0$ in the supermanifold $V$ and vanishing at $0$.
A morphism $U\colon(V,Q)\to (V',Q')$ of $Q$\ndash manifolds is a coalgebra
morphism $C^0(V)\to C^0(V')$ of degree $0$
obeying $Q'\circ U=U\circ Q$.

In explicit terms, a coderivation $Q$ of $C^0(V)$ is uniquely determined by
its composition $p_1\circ Q$ with the canonical projection $p_1\colon C(V)\to V$ 
sending $I_j(V)$ to $0$ for $j\neq 1$, see Lemma \ref{l-Alighieri}. The restriction
of $p_1\circ Q$ to $I_j(V)=(V^{\otimes j})^{S_j}$ is a map
$Q_j\colon I_j(V)\to V$ of degree 1, the $j$th Taylor component. 
The
condition $[Q,Q]=0$ is then equivalent to 
\[\sum_{j+k=n}\sum_{l=0}^{j-1}
 Q_j\circ(1^{\otimes l}\otimes Q_k\otimes 1^{\otimes( j-l-1)})=0,\]   
on $I_n(V)$, $n=1,2,\dots$.
Similarly, a coalgebra morphism $U$ is uniquely determined by its
Taylor components $U_j=p_1\circ U|_{I_j(V)})$, $(j\geq 1)$. The
$Q$\ndash manifold morphism property is then
\begin{eqnarray*}
\lefteqn{\sum_{j_1+\cdots+j_k=n}Q_k\circ(U_{j_1}\otimes\cdots\otimes U_{j_k})}
\\
&&= \sum_{j+k-1=n}\sum_{l=0}^{j-1}
 U_j\circ(1^{\otimes l}\otimes Q_k\otimes 1^{\otimes(j-l-1)})
\end{eqnarray*}
on $I_n(V)$.

If $\mathfrak g$ is a differential
graded Lie algebra, then $\mathfrak g[1]$ (with $\mathfrak g[1]^i=\mathfrak g^{i+1}$) is a $Q$\ndash manifold: the Taylor components of the coderivation vanish
except $Q_1$ and $Q_2$, which are given in terms of the differential
$d$ and the bracket by
\[
Q_1=d,\qquad Q_2(\gamma_1,\gamma_2)=(-1)^{|\gamma_1|}
[\gamma_1,\gamma_2],\qquad \gamma_i\in\mathfrak{g}[1]^{|\gamma_i|}=\mathfrak g^{|\gamma_i|+1},
\]
(a more pedantically correct notation for the right\ndash hand side of the equation for $Q_2$ would be
 $(-1)^{|\gamma_1|}s^{-1}[s\gamma_1,s\gamma_2]$).

\medskip

\noindent{\bf Definition.}
A flat $L_\infty$\ndash algebra structure on a vector space
(or $R$\ndash module) $\mathfrak g$ is a $Q$ manifold 
structure on $\mathfrak g[1]$.

\medskip

It is convenient to express a flat $L_\infty$\ndash algebra
in terms of 
the structure maps $\tilde Q_j\in\mathrm{Hom}^{2-j}
(\wedge^j\mathfrak g, \mathfrak g)$ (differential and
higher Lie brackets)
of $\mathfrak g$. They are the $Q_j$ up to sign. 
The precise relation is
most easily expressed using the tautological map $s\colon\mathfrak g[1]\to\mathfrak
g$ of degree $1$. We then have $\tilde Q_j=s^{-1}\circ Q_j\circ(s\otimes\cdots
\otimes s)$. Note that $s^{\otimes j}$ intertwines the action of 
$S_n$ on $\mathfrak{g}[1]^{\otimes j}$ with the alternating action
of $S_n$ on $\mathfrak{g}^{\otimes j}$:
$s^{\otimes j}\sigma=\mathrm{sign}(\sigma) \sigma s^{\otimes j}$, $s\in S_j$.
Thus if $Q_j$ are symmetric $\tilde Q_j$ are skew\ndash symmetric.
Explicitly,
\[
Q_j(\gamma_1,\dots,\gamma_j)=(-1)^{\sum_{i=1}^j(j-i)|\gamma_i|}
s^{-1}\tilde Q_j(s\gamma_1,\dots,s\gamma_j),\qquad \gamma_i\in\mathfrak g[1].
\]

\subsection{The local formality theorem}
\begin{thm}\label{t-form}
 Let $A=\mathbb R[x_1,\dots,x_d]$ be the graded commutative algebra
generated by $x_i$ of degree $\eps_i$, $i=1,\dots,d$ or its completion
$A=\mathbb R[[x_1,\dots,x_d]]$. Then there exists a morphism of 
formal pointed $Q$\ndash manifolds
$U\colon\mathcal{T}(A)[1]\to\mathcal{D}(A)[1]$
 such that $U_1$ is the HKR quasiisomorphism.
\end{thm}

\subsection{The construction of the $L_\infty$\ndash morphism}
The condition for the Taylor components $(U_n)_{n\geq1}$ of a
morphism of $Q$-manifolds 
$U$ from $\mathcal{T}(A)[1]\to \mathcal{D}(A)[1]$ are simplified
if we add a component $U_0=\mu_A$, the product in $A$. Geometrically,
this means that we shift the origin of the $Q$-manifold $\mathcal{D}(A)$ by $U_0$
to a point where the vector field $Q$ is purely quadratic.
The conditions for $(U_n)_{n\geq 0}$ 
to be satisfied are then
\begin{equation}\label{e-formality}
\sum_{n_1+n_2=n}Q_2\circ(U_{n_1}\otimes U_{n_2})=
\sum_{l=0}^{n-2}
 U_{n-1}\circ(1^{\otimes l}\otimes Q_2\otimes 1^{\otimes(n-l-2)}),\qquad n\geq 1,
\end{equation}
on the space of symmetric tensors $I_n(\mathcal{T}(A)[1])$.
Replacing $\oplus I_n$ by the isomorphic shuffle coalgebra (see \ref{s-free}),
 we
can write this as
\begin{eqnarray}\label{e-Ticknor}
\lefteqn{\sum_{p+q=n}\sum_{(p,q)-\mathrm{shuffles}}
\pm Q_2(U_{p}(\gamma_{\sigma(1)},\dots,
\gamma_{\sigma(p)}),U_{q}(\gamma_{\sigma(q+1)},\dots,
\gamma_{\sigma(n)}))}\notag
\\
&&=
 \sum_{i<j}(-1)^{\eps_{ij}}U_{n-1}( Q_2(\gamma_i,\gamma_j),\gamma_1,\dots,\hat\gamma_i,\dots,
\hat\gamma_j,\dots,\gamma_n),\\
\eps_{ij}&=&\sum_{k=1}^{i-1}|\gamma_i|\cdot|\gamma_k|+\sum_{k=1}^{j-1}|\gamma_j|\cdot|\gamma_k|-|\gamma_i|\cdot|\gamma_j|,\notag
\end{eqnarray}
and $\pm$ is the same sign as in \eqref{e-sign}.
The maps $U_n$ are  sums of integrals over configuration
spaces on the upper half\ndash plane $H$. 
Let $n,m$ be non-negative integers such that
$2n+m\geq 2$. Let $\mathrm{Conf}^+_{n,m}=\{(z,x)\in H_+^n\times \mathbb R^m\,|\, z_i\neq z_j\, (i\neq j),
x_1<\dots<x_m\}$, with orientation form $d^2z_1\cdots d^2z_n
 dx_1\cdots dx_m$, where $d^2z_i= d\,\mathrm{Re}(z_i)\,
d\,\mathrm{Im}(z_i)$. The group $G_2$ of affine transformations 
$z\to \lambda z+a$, $\lambda>0,a\in\mathbb R$ acts freely
on $\mathrm{Conf}^+_{n,m}$
since $2n+m\geq 2$ and preserves the orientation. Let the orientation of $G_2$
be defined by the volume form $da\wedge d\lambda$.  
The quotient $\mathcal{C}^+_{n,m}$ is
oriented as in \cite{AMM} in such a way that any
 trivialisation $G_2\times \mathcal{C}^+_{n,m}\to \mathrm{Conf}^+_{n,m}$ of the left  principal $G_2$\ndash bundle $\mathrm{Conf}_{n,m}$ is orientation preserving. Here is an explicit description:
if $n\geq1$,
$\mathcal{C}^+_{n,m}$ may be identified with the submanifold
of $\mathrm{Conf}^+_{n,m}$ consisting of points with $z_1=i$
and orientation form $d^2z_2\cdots d^2z_n dx_1\cdots dx_m$. If $m\geq2$ it
can be identified with the submanifold given by $x_1=0$,
$x_m=1$, with orientation form $(-1)^m d^2z_1\cdots d^2z_n
dx_2\cdots d x_{m-1}$. Let $\mathcal{G}_{n,m}$ be the
set of graphs $\Gamma=(V_\Gamma,E_\Gamma)$ with the following
properties: the set of vertices 
$V_{\Gamma}=\{1,\dots,n,\bar 1,\dots,\bar m\}$ consists of
vertices of the first type $1,\dots,n$ and vertices of
the second type $\bar1,\dots,\bar m$; the edges $(i,j)\in
E_\Gamma\subset V_\Gamma\times V_\Gamma$  are such that
$i$ is always of the first type and $i\neq j$.

Let moreover $\tau\in\mathrm{End}(\tilde A\otimes\tilde A)$ be
the endomorphism
\[
\tau =\sum_{\alpha=1}^{d}(-1)^{\eps_\alpha}
\partial_{\theta_\alpha}\otimes \partial_{x_\alpha}.
\]
To each graph $\Gamma\in\mathcal{G}_{n,m}$ we associate an element
$\omega_\Gamma$
of the graded algebra $\Omega(\mathcal{C}^+_{n,m})\otimes\mathrm{End}(A_{n,m})$
of  differential forms with values in the endomorphisms of
\[
A_{n,m}=\tilde A^{\otimes n}\otimes A^{\otimes m}.
\] 
Let the factors in the tensor product $A_{n,m}$
be numbered $1,\dots,n,\bar1,\dots,\bar m$ and  for $i,j\in V_\Gamma$ let 
$\tau_{ij}\in\mathrm{End}(A_{n,m})$  the endomorphism acting as $\tau$
on the factors $i$ and $j$  and as the identity on the other factors:
\[
\tau_{ij}=\sum_{\alpha}(-1)^{\eps_{\alpha}}1\otimes\cdots\otimes 1
\otimes
\partial_{\theta_\alpha}\otimes1\otimes\cdots\otimes 1\otimes\partial_{x_\alpha}
\otimes
1\otimes\cdots\otimes 1.
\]
Set $d\phi_{ij}\in\Omega^1(\mathcal{C}^+_{n,m})$ the differential of
the Kontsevich angle function 
$\phi(z_i,z_j)=\arg(z_i-z_j)-\arg(\bar z_i-z_j)$ (with $z_{\bar k}=x_k$ if
$\bar k$ is of the second type).
Both $\tau_{ij}$ and $d\phi_{ij}$ are elements of the algebra
$\Omega(\mathcal{C}^+_{n,m})\otimes\mathrm{End}(A_{n,m})$ of degree 1 and
$-1$, respectively. Thus their product is of degree $0$ and
\[
\omega_\Gamma=\prod_{(i,j)\in E_\Gamma}d\phi_{ij}\tau_{ij}
\]
is independent of the choice of ordering of factors.

Let $U_n$ be the map 
$ \oplus_{m\geq 0}(\tilde A^{\otimes n}\otimes A^{\otimes m})\to A$ 
\begin{equation}\label{e-Un}
U_n=\sum_{m\geq0}
(-1)^{(\sum|\gamma_i|-1)m}\sum_{\Gamma\in\mathcal{G}_{n,m}}U_\Gamma,
\end{equation}
where
\[
U_\Gamma=\mu
\int_{\mathcal{C}^+_{n,m}}\prod_{(i,j)\in E_\Gamma}d\phi_{ij}\tau_{ij}
\]
is the composition
\[
(\tilde A^{\otimes n}\otimes A^{\otimes m})\stackrel{\omega_\Gamma}{\rightarrow} 
\Omega(\mathcal{C}^+_{n,m})\otimes(\tilde A^{\otimes n}\otimes A^{\otimes m})\to
\tilde A^{\otimes n}\otimes A^{\otimes m}\stackrel{\epsilon\mu}{\to} A
\]
The second map is the integration $\omega\otimes a\to (\int\omega) a$ 
(of degree $-2n-m+2$) and
is defined to be zero on differential forms of the wrong degree.
Thus $U_\Gamma$ on $\tilde A^{\otimes n}\otimes A^{\otimes m}$
vanishes unless the number of edges is
\[
|E_\Gamma|=2n+m-2.
\]
The map $\epsilon\mu\colon \tilde A^{\otimes n}\otimes A^{\otimes m}\to A$ is the product
in $\tilde A$ followed by the projection $\epsilon\colon\tilde A\to A$ sending $\theta_i$ to $0$.

\begin{proposition}\label{t-Li} The maps $U_n$ are Taylor components of a morphism 
$U$ with the properties stated in Theorem \ref{t-form}.
\end{proposition}

The proof of this theorem is based on the Stokes theorem as in \cite{Kontsevich}.
The quadratic relation \eqref{e-formality} are obtained from a sequence of
relations for integrals over configuration spaces associated to graphs. Let
$\Gamma\in G_{n,m}$ be a graph such that 
\[
|E_\Gamma|=2n+m-3.
\]
Then $d\omega_\Gamma$ (which vanishes) is a form of degree $2n+m-2$ and we
have the Stokes theorem on the Kontsevich compactification $\bar C^+_{n,m}$
of $C^+_{n,m}$:
\[
0=\int_{\bar C^+_{n,m}}d\omega_\Gamma=\sum_{i}\int_{\partial_i\bar C^+_{n,m}}\omega_\Gamma.
\]
The sum is over the faces of the manifold with corners $\bar C^+_{n,m}$. The
faces contributing non-trivially are of two types. 

(a) Faces of the first type
are diffeomorphic to $C^+_{n',m'}\times C^+_{n'',m''}$ with $n'+n''=n$ and
$m'+m''=m+1$ and correspond to limiting configurations where $n'$ points in $H$
and $m'$ consecutive points with labels $\overline{l+1},\dots,\overline{l+m'}$
on the real line converge to a single point.
The orientation from the Stokes theorem differs from the product orientation
by a factor $(-1)^{lm'+l+m'}$, as computed in \cite{AMM}. (b) Faces of the
second type are diffeomorphic to $C_{n'}\times C^+_{n'',m}$. with $n'+n''=n-1$  
and correspond to limiting configurations where $n'$ points in $H$ converge
to the same point in the upper half plane. Here the relative position of
these $n'$ collapsing points is parametrised by the manifold $C_{n'}$ ($n'\geq 2$), 
the quotient
of $\mathbb C^{n'}\setminus
\cup_{i<j}\{z_i=z_j\}$ by the group $G_3$ of
affine transformations $z\mapsto \lambda z+a$ with $\lambda>0$ and $a\in\mathbb C$
and orientation form $d^2a\wedge d\lambda$. 
By Kontsevich's lemma (see \cite{Kontsevich}, Lemma 6.6)
the integrals over these faces vanish except
for $n'=2$.
In this case the induced orientation
on the face differs by the product orientation by a factor $-1$, see \cite{AMM}. 

The faces of the first type will contribute to the expression on the left-hand side
of \eqref{e-Ticknor} as in \cite{Kontsevich}. We need to keep track of the signs.
Let us consider the case of a face of the first type
in which the $n'$ points in the upper half-plane collapsing to a point on the
real axis are the last $n'$ points; this corresponds to the trivial shuffle in
\eqref{e-Ticknor}. The remaining shuffles are treated similarly or can be related
to the trivial one by permutation symmetry considerations. Let us denote
accordingly $|\gamma''|=\sum_{i=1}^{n''}|\gamma_i|$, $|\gamma'|=\sum_{i=n''+1}
^{n}|\gamma_i|$, $|\gamma|=|\gamma'|+|\gamma''|=\sum_1^n|\gamma_i|$.
The sign with which the integral of $\omega_\Gamma$ over this face
 contributes to the term
\begin{equation}\label{e-Mannings}
U_{n''}(\gamma_1,\dots,\gamma_{n''})(1^{\otimes l}\otimes U_{n'}(\gamma_{n''+1},\dots,\gamma_n)\otimes 1^{\otimes{m''-1}}),
\end{equation}
appearing (with a certain sign we give below) in the left-hand
side of \eqref{e-Ticknor} is 
\[
(-1)^{lm'+m'+l}(-1)^{|\gamma''|m'}(-1)^{-(|\gamma'|-1)m'-(|\gamma''|-1)m''}
\]
The first sign of this product
 comes from comparing the orientations, as discussed above,
the second from moving $\gamma_i$, $i\leq n''$
to the left of $\prod_{ij\in E_{\Gamma'}}\tau_{ij}$ ($|E_{\Gamma'}|\equiv m'
\mod 2$), the third appears in the definition of $U_{n'}$, $U_{n''}$.
The same term \eqref{e-Mannings}
appears in the left-hand side of \eqref{e-Ticknor} with a
sign 
\[
(-1)^{|\gamma''|-1}(-1)^{(|\gamma'|-1)(m''-1)+(m'-1)l},
\]
which is the product of the sign coming from comparing $Q_2$ to the Gerstenhaber
bracket and the sign appearing in the definition of the Gerstenhaber bracket. The
ratio between these signs is
\[
(-1)^{|\gamma|(m'+m''-1)}=(-1)^{|\gamma|m},
\]
which is the sign with which the considered face contributes to the left-hand side
of \eqref{e-Ticknor}. Let us turn to the right-hand side: the face in which
the first two points in the upper half-plane collapse contributes to the
term right-hand side of \eqref{e-Ticknor} with $i=1,j=2$ with the same sign
$(-1)^{|\gamma|m}$ which is the sign appearing in the definition of $U_{n-1}$
by taking into account the fact that $|[\gamma_1,\gamma_2]|+\sum_{i\geq 3}|\gamma_i|
=|\gamma|-1$. The orientation sign $(-1)$ is used to write the term on the right-hand
side, and no other sign appears since the expression
\[
\mu\circ\tau(\gamma_1\otimes\gamma_2)=(-1)^{|\gamma_1|-1}
\sum_i \gamma_1\!
\stackrel{\leftarrow}\partial_{\theta_i}
\stackrel{\rightarrow}\partial_{x_i}\!\!
\gamma_2
\]
obtained from Stokes has the same sign as the corresponding term in
$Q_2(\gamma_1,\gamma_2)=(-1)^{|\gamma_1|-1}[\gamma_1,\gamma_2]$.

There remains to show that $U_1$ is the HKR map. Let
\[
\gamma=\gamma^{i_1\cdots i_m}\theta_{i_1}\cdots\theta_{i_m}.
\]
Then there is only one graph in $\mathcal{G}_{1,m}$ contributing to $U_1(\gamma)$.
Its edges are $(1,\bar j)$, $1\leq j\leq m$. We then have for any $f\in A^{\otimes m}$.
\begin{eqnarray*}
U_1(\gamma)f
&=&(-1)^{(|\gamma|-1)m}
\int_{C^+_{1,m}}\prod_{i=1}^m 
\omega_{1,\bar i}
(\gamma\otimes f)\\
&=& (-1)^{(|\gamma|-1)m}\frac{1}{(2\pi)^m}
\int_{0<\varphi_1<\cdots<\varphi_m<2\pi}
\prod_{i=1}^m 
d\varphi_i\tau_{1\bar i}
(\gamma\otimes f)\\
&=&(-1)^{(|\gamma|-1)m+m(m-1)/2}\frac1{m!}
\tau_{1\bar 1}\cdots\tau_{1\bar m}
(\gamma\otimes f)\\
&=& 
(-1)^{(|\gamma|-1)m+m(m-1)/2}
(-1)^{|\gamma|m+\sum_{\alpha=1}^m(\alpha(1-\epsilon_{i_\alpha})+\epsilon_{i_\alpha})}
\\
&&\gamma^{i_1\cdots i_m}
\mu\circ(\partial_{i_1}\otimes \cdots\otimes\partial_{i_m})
f\\
&=&(-1)^{\sum_{\alpha=1}^m(\alpha-1)\epsilon_{i_\alpha}}\gamma^{i_1\cdots i_m}
\mu\circ(\partial_{i_1}\otimes \cdots\otimes\partial_{i_m})
f
\end{eqnarray*}

\subsection{Globalisation}
The $L_\infty$\ndash morphism of Proposition \ref{t-Li} obeys all the additional
properties of \cite{Kontsevich} 
needed to go from a local to a global $L_\infty$\ndash morphism. In particular
the non-trivial
fact that $U_n$ vanishes if one of its arguments is a linear vector field
is valid here because of the vanishing of the same integrals over configuration
spaces as in \cite{Kontsevich}. One can then deduce the existence of a
morphism of $Q$\ndash manifolds $U\colon\mathcal{T}(A)[1]\to \mathcal{D}(A)[1]$ for the
algebra $A$ of functions on any supermanifold along the lines of \cite{Kontsevich}, \cite{Kontsevich2}
or \cite{Dolgushev}.

\newcommand{\calT}{{\mathcal{T}}}
\newcommand{\calD}{{\mathcal{D}}}
\newcommand{\de}{\partial}
\newcommand{\sfA}{{\mathsf{A}}}
\newcommand{\sfB}{{\mathsf{B}}}
\newcommand{\bU}{{\mathbf{U}}}

The explicit contruction goes as follows. To fix notations, let $E\to M$ be the graded vector bundle that
realizes the given supermanifold; viz., $A=\Gamma(M, \Hat SE^*)$. Using a connection, we may
identify $\mathcal{T}(A)$ with $\Gamma(M,ST[-1]M\Hat\otimes\Hat SE^*\Hat\otimes\Hat SE[-1])$ as GLAs.

We may get a Fedosov resolution thereof following \cite{Dolgushev}. A simplification is obtained by using an idea contained
in \cite{Bordemann1}. Namely, we consider the complex $\Omega^\bullet(\calT):=\Gamma(\Lambda^\bullet T^*M\otimes\calT)$,
with 
\[
\calT=S(T^*M)\Hat\otimes ST[-1]M\Hat\otimes\Hat SE^*\Hat\otimes\Hat SE[-1]
\]
as a GLA. This is the only difference
with the construction of \cite{Dolgushev}. The rest goes exactly the same way. Namely,
one constructs a compatible differential $D$ with cohomology concentrated in degree zero and equal to $\mathcal{T}(A)$
as follows. First one considers the (globally well-defined)
differential $\delta:=[d x^i\frac\de{\de y^i},\ ]$, with $\{x^i\}$ local coordinates on $M$ and $\{y^i\}$
the induced local coordinates on the tangent fibers. Then one picks up a torsion-free affine connection which together with
the already chosen connection on $E$ defines a connection on $\calT$. The induced covariant derivative $\nabla$ commutes with $\delta$
since the connection is torsion free. One then kills the curvature
constructing by induction, exactly as in \cite{Dolgushev}, an element $\sfA$ of
$\Omega^1(\calT)$ such that $D:=\nabla-\delta+A$ squares to zero and has the wished-for properties.
Similarly one gets a Fedosov resolution $\Omega^\bullet(\calD)$ of $\mathcal{D}(A)$.

Next, on a coordinate neighbourhood $W$, one defines  a splitting $D=d+\sfB$ with $d=d x^i\frac\de{\de x^i}$.  
The local $L_\infty$-quasiisomorphism defined in the previous subsection,
may be extended to an $L_\infty$-quasiisomorphism 
\[
U_W\colon (\Omega^\bullet(\calT)|_W,d,[\ ,\ ])\leadsto(\Omega^\bullet(\calD)|_W,d+\de,[\ ,\ ]),
\]
with $\de$ the Hochschild differential. Now one observes that
$\sfB$ is a MC element in $\Omega^\bullet(\calT)|_W$ (as a vector field) and $\Omega^\bullet(\calD)|_W$ (as a first-order differential operator) and
that it it is mapped to itself by $U_W$ since it is a vector field. 
So one can localize $U_W$ at $\sfB$ and get a new $L_\infty$-quasiisomorphism
\[
\bU_W\colon (\Omega^\bullet(\calT)|_W,D,[\ ,\ ])\leadsto(\Omega^\bullet(\calD)|_W,D+\de,[\ ,\ ]).
\]
Since $\sfB$ transforms from one coordinate neighbourhood to another by the addition of a linear vector field and since higher components of the
$L_\infty$-quasiisomorphism vanish on linear vector fields, one realizes that $\bU_W$ does not really depend on a choice of local
coordinates, so it extends to a global $L_\infty$-quasiisomorphism $\bU$. Finally, one may modify $\bU$ in such a way that
its image lies in the DGLA of zero $D$-cochains, which may be identified with $\mathcal{D}(A)$.


\begin{thebibliography}{BGHHW}
\bibitem[AMM]{AMM} D. Arnal, D. Manchon, M. Masmoudi, {\em Choix des signes pour la formalit\'e de M. Kontsevich.}  Pacific J. Math.  203  (2002),  no. 1, 23--66. 
\bibitem[B1]{Bordemann1} 
M. Bordemann, {\em On the Deformation Quantization of Super-Poisson Brackets.}
\texttt{q-alg/9605038}.
\bibitem[B2]{Bordemann} 
M. Bordemann, {\em (Bi)modules, morphismes et r\'eduction des star-produits~: le cas symplectique, feuilletages et obstructions.}
\texttt{math.QA/0403334}.
\bibitem[BGHHW]{BGHHW}  M. Bordemann, G. Ginot, G. Halbout, H.-C. Herbig, S. Waldmann,
{\em
Star-repr\'esentations sur des sous-variet\'es co-isotropes.} \texttt{math.QA/0309321}.
\bibitem[CF]{CaFe2003}
A. S. Cattaneo, G. Felder,
{\it Coisotropic submanifolds in Poisson geometry and branes in the Poisson sigma model.} \texttt{math.QA/0309180}.
\bibitem[D]{Dolgushev}
V. Dolgushev, {\em Covariant and equivariant
formality theorems.} Adv. Math. {\bf 191} (2005), no. 1, 147--177.
\bibitem[Ge]{Gerstenhaber}
M. Gerstenhaber,
{\em The cohomology structure of an associative ring.}
Ann. Math. {\bf 78} (1963), no. 2, 267--288.
\bibitem[HKR]{HKR}
G. Hochschild, B. Kostant, A. Rosenberg,
{\em Differential forms on regular affine algebras.}
Trans. Amer. Math. Soc. 102 (1962) 383--408.
\bibitem[K1]{Kontsevich}
M. Kontsevich,
{\em Deformation quantization of Poisson manifolds.}  Lett. Math. Phys.  {\bf 66}  (2003),  no. 3, 157--216.
\bibitem[K2]{Kontsevich2}
M. Kontsevich,
{\em Deformation quantization of algebraic varieties.}  Lett. Math. Phys.  {\bf 56}  (2001),  no. 3, 271--294.
\bibitem[OP]{OP} Y.-G. Oh, J.-S. Park, 
{\em Deformations of coisotropic submanifolds
and strongly homotopy Lie algebroid.}
 \texttt{math.SG/0305292}.
\bibitem[LS]{LS}
 S.L. Lyakhovich, A.A. Sharapov,
{\em BRST theory without Hamiltonian and Lagrangian.} \texttt{hep-th/0411247}.
\bibitem[MK]{Mackenzie}
K. Mackenzie,
{\em Lie groupoids and Lie algebroids in differential geometry.}
London Mathematical Society Lecture Note Series, 124.
Cambridge University Press, Cambridge, 1987.
\bibitem[M]{MacLane}
S. Mac Lane, {\em Homology}, third edition, Springer 1975
\bibitem[R]{Roytenberg}
D. Roytenberg, {\em
 Courant algebroids, derived brackets and even symplectic supermanifolds.} PhD thesis, UC Berkeley, 1999, \texttt{math.DG/9910078}.
\bibitem[S1]{Stasheff63}
J. Stasheff,
{\em Homotopy associativity of H-spaces I,II.} Trans. Amer. Math. Soc. {\bf 108} (1963), 275--312.
\bibitem[S2]{Stasheff93}
J. Stasheff,
{\em The intrinsic bracket on the deformation complex of an associative algebra.} J. Pure Appl. Math. {\bf 89} (1993), 231--235.
\bibitem[Vo]{Voronov}
Th. Voronov,
{\em Higher derived brackets and homotopy algebras.}
\texttt{math.QA/0304038}.
\bibitem[Ve]{Vey} J. Vey,
{\em D\'eformation du crochet de Poisson sur
une vari\'et\'e symplectique.} Comm. Math. Helv.
{\bf 50} (1975), 421--454.
\bibitem[W]{Weinstein}
A. Weinstein,
{\em Coisotropic calculus and Poisson groupoids.}
J. Math. Soc. Japan {\bf 40} (1988), 705--727.
\end{thebibliography}
\end{document}